\documentclass[preprint, 12pt, 3p, authoryear]{elsarticle}
\newdimen\Lmargin
\newdimen\Rmargin
\Rmargin=\paperwidth \advance\Rmargin by-\Lmargin \advance\Rmargin by-\textwidth
\ifdim\Lmargin>\Rmargin \Rmargin=\Lmargin \fi
\def\Scentering{\leftskip=0pt plus 1fil minus \Rmargin
                \rightskip=\leftskip}

\usepackage{amssymb}
\usepackage{graphicx}
\usepackage[colorlinks=true]{hyperref}
\usepackage{multirow}
\usepackage{booktabs}
\usepackage{natbib}
\usepackage{tabularx}
\usepackage{booktabs}
\usepackage{makecell}
\usepackage{amsmath}
\usepackage{lscape}
\usepackage{utfsym}
\usepackage{longtable}
\usepackage{caption}
\usepackage{amsthm}

\begin{document}

\begin{frontmatter}
\newtheorem{theorem}{Theorem}
\newtheorem{proposition}{Proposition}
\newtheorem{corollary}{Corollary}
\newtheorem{lemma}{Lemma}
\newdefinition{example}{Example}
\newdefinition{definition}{Definition}
\newdefinition{remark}{Remark}


\title{Generalized Ordinal Priority Approach for Multi-Attribute Decision-Making under Incomplete Preference Information}

\author[inst1]{Renlong Wang\corref{cor1}}
\ead{13127073530@163.com}
\affiliation[inst1]{organization={School of Emergency Management Science and Engineering, University of Chinese Academy of Sciences},
            city={Beijing},
            postcode={100049}, 
            country={China}}
            
\cortext[cor1]{Corresponding author}

\begin{abstract}

The Ordinal Priority Approach (OPA) is a multi-attribute decision-making (MADM) method to determine the relative importance (weights) of experts, attributes, and alternatives. This study formally establishes the fundamental properties of OPA, including solution efficiency, analytical solution expression, the decomposability of optimal decision weights, and its relationship with rank-based surrogate weights. Building on these properties, we propose a Generalized Ordinal Priority Approach (GOPA) based on an “estimate-then-optimize” contextual optimization framework for MADM when preference information is incomplete. In the first stage, we derive utility distributions for ranked alternatives in discrete and continuous prospects by minimizing cross-entropy utility under partial preference information, including weak order relations, absolute differences, ratio scales, and lower bounds. Rank-based surrogate weights and risk preference utility functions serve as the global utility structure for discrete and continuous prospects, respectively. The elicited utility information is then introduced into the second-stage problem to simultaneously optimize the weights of experts, attributes, and alternatives within a normalized weight space. Metrics for validating the group decision outcomes of GOPA, including percentage standard deviation, correlation coefficient, and confidence level measurement, are proposed. Theoretical analysis reveals several advantageous properties of GOPA, including model generalizability, analytical solvability, and risk preference independence. Furthermore, this study provides a lower bound reference for transforming the general optimization-based weight elicitation problems into optimization problems with stochastic dominance constraints. The applicability of GOPA is demonstrated through an improvisational emergency supplier selection problem.

\end{abstract}



\end{frontmatter}

\section{Introduction}
\label{sec-1}

In the past decades, multi-attribute decision-making (MADM) has been regarded as an effective tool for addressing complex decision problems, adept at tackling challenges such as objective conflicts, data diversity and high uncertainty \citep{R.L.C23, W24}. A classical problem involves selecting an optimal alternative or obtaining a global ranking from a set of alternatives under multi-attributes based on multi-expert opinions. More specifically, for a given expert set $\mathcal{Q}$, attribute set $\mathcal{N}$, and alternative set $\mathcal{M}$, the evaluation score $V_{k}$ of alternative $k$ can be determined by a mapping $F: \mathbb{R}^{I} \times \mathbb{R}^{J} \mapsto \mathbb{R}$ that has an associated collection of the weights for experts $w_{i}$ and attributes $w_{j}$:
\[
Z_k=F( v_{ijk} ) = \sum_{j\in\mathcal{N}} w_j \sum_{i\in\mathcal{Q}} w_i v_{ijk}, \quad \forall k \in \mathcal{M},
\]
where $v_{ijk}$ is the performance score of alternative $k$ under attribute $j$ given by expert $i$. Typically, much of the research is based on precise weights of experts and attributes obtained through sophisticated heuristic methods \citep{Z.K.Y.Z.D23}. Meanwhile, some studies have emerged to elicit weight through incomplete information, which is often referred to as imprecise or partial preference information, to address challenges such as time constraints, inadequate data and domain knowledge, and limitations in decision-makers' attention and information processing capabilities \citep{A24}. Of course, if decision-makers can provide all the information required to solve MADM problems, prior (sophisticated) methods based on precise data is recommended. 

Current methods for weight elicitation under incomplete preference information mainly fall into two categories: optimization-based methods and extreme point-based methods. The former aims to model and solve the optimal weight combination under the constraints of incomplete preference information \citep{K.A99, M00, W.C.Y07}. Conversely, the latter capitalizes on extreme points within the incomplete data to ascertain weights \citep{A15, A17, A24}. Some research have suggested an evolving synthesis of these two approaches. Within optimization-based approaches, the DEA preference voting model, which incorporates partial preference information in a ratio scale, is prevalently employed \citep{W.C.Y07}. Correspondingly, \citet{A17} has derived a closed-form expression for the optimal solution of the DEA preference voting model by imposing weight constraints from extreme points of partial preference information. Furthermore, \citet{A24} has derived a dual problem of linear programming to obtain closed-form solutions and establish a strict ranking of extremal points of attribute weights. Additionally, a prevalent weight elicitation approach involves rank-based surrogate weights for weight elicitation, where each rank-based surrogate weight is uniquely represented by a set of extreme points \citep{B.N23}. Thus, \citet{L24} propose explicit expressions for weights of various simplex centroids in ranking voting frameworks inspired by specific simplex centroids of ROC weights. However, previous studies considering incomplete preference information primarily focused on determining attribute or expert weights, neglecting simultaneous consideration of both. This necessitates additional techniques for integrating expert opinions to facilitate group decision-making. Meanwhile, evaluating the performance scores of alternatives across various attributes by experts is also a challenging task in subjective decision-making. Moreover, these studies do not account for the influence of decision-makers' risk preferences on weight assignment, often assuming all decision-makers have similar risk-neutral characteristics. It is critical to recognize the potential impact of evaluation scale distortions resulting from differing risk preferences on decision outcomes \citep{T.K74}.

The Ordinal Priority Approach (OPA) is an emerging optimization-based method for addressing MADM under incomplete information \citep{A.M.F.L20}. OPA derives weights from a normalized weight space with ordinal preference using linear programming, simultaneously determining weights for experts, attributes, and alternatives. OPA is applicable to both group and individual decision-making. Compared to typical inputs in other MADM methods, such as cardinal values and pairwise comparisons, ordinal data is more accessible to obtain and more stable. Unlike traditional methods, OPA avoids the need for data standardization, expert opinion aggregation, and prior weight acquisition. Additionally, various extensions of OPA have emerged, such as fuzzy OPA \citep{S.M.D.L23, Z.H.L24}, rough OPA \citep{D.L.J.G.C23, K.P.D.E.D23}, Grey OPA \citep{M.D.J.Z21}, and robust OPA \citep{M.A.D22} for handling decision data uncertainty, DGRA-OPA \citep{W24} and TOPSIS-OPA \citep{M.A.Y.L22} for addressing large-scale group decision-making, and partial OPA \citep{W.S.C.S.C.G24} for managing potential Pareto dominance. Consequently, OPA has garnered increasing attention recently and has been applied across various domains, including sustainable transportation evaluation \citep{P.D.G.T.K22}, supplier evaluation \citep{M.D.J.Z21, M.J22}, blockchain analysis \citep{ S.M.D.L23, Z.H.L24}, project portfolio selection\citep{M.A.D22}, and emergency recovery planning \citep{W24}.  However, as a practical approach, research on the fundamental properties of OPA remains scarce. The utilization of ordinal data can only capture dominance relations within decision-makers’ preferences, failing to consider the partial preference information in the form such as ratio scale and absolute difference. Consequently, OPA cannot characterize the different weight distributions formed by the risk preferences of decision-makers, which are crucial for real-world decision-making. 

Overall, MADM with incomplete information poses the following challenges: i) insufficient research on the foundational theory of OPA; ii) inability to consider the impact of risk preferences and their heterogeneity on decisions; and iii) the necessity of pre-obtaining weights of experts and attributes. Therefore, we first formally establish the fundamental mathematical properties of OPA. Based on these properties, we introduce a Generalized Ordinal Priority Approach (GOPA) using a two-stage “estimate-then-optimize” contextual optimization framework to address decision weight elicitation challenges in MADM with incomplete information (covering experts, attributes, and alternatives). In the first stage, we introduce utility functions derived from partial preference information and risk attitude, examining the application of the cross-entropy utility minimization method to handle preference information in limited data under discrete and continuous prospects. The utility distribution derived in the first stage serves as uncertain parameters for optimizing decision weights in the second stage. We analyze the properties of GOPA, including the degradation of cross-entropy utility maximization, risk preference within the contextual optimization framework, and the analytical solution of GOPA. We develop metrics to validate GOPA group decision results. Finally, we discuss the formulation of GOPA to convert optimization-based utility elicitation methods into a classical optimization problem with stochastic dominance constraints and propose formulas for correcting elicitation errors.

The primary contributions of this study can be summarized as follows:
\begin{itemize}
    \item [1)] This study formally establishes several fundamental properties of OPA, including the solution efficiency, analytical solution and its decomposability, and the relationship between optimal weights and rank-based surrogate weights. The decomposability provides a theoretical foundation for setting reference points in the prospect theory-based extension of OPA. These findings could deepen the understanding and provide a foundation for advancing theoretical research of OPA.
    \item [2)] This study proposes GOPA, extending the rank order centroid weights within OPA formulation to a more general utility form, to tackle the challenge in MADM with incomplete preference information. GOPA employs an “estimate-then-optimize” contextual optimization framework that integrates general utility elicitation across discrete and continuous prospects, encompassing global utility structures with risk attitude and partial preference information. This accommodates the personalized preferences of experts across various attributes. GOPA demonstrates advantageous properties in the decision science domain, including model generalizability, analytical solvability, and risk preference independence.
    \item [3)] GOPA provides a general two-stage framework for optimization-based MADM methods considering the risk preference of decision-makers: eliciting the utility of alternatives through risk preferences and partial preference information across attributes, followed by optimizing weights within a normalized weight space encompassing experts, attributes, and alternatives. Furthermore, the decision weight expression in the analytical solution of GOPA can serve as a universal utility expression for other utility elicitation research by inputting the personalized utility distributions of alternatives.
    \item[4)] The objective function value in the analytical solution of GOPA provides a lower bound reference guided by the utility from subjective decision information, representing a novel exploration in this domain. This forms the foundation for converting subjective preference-based utility elicitation into optimization with stochastic dominance constraints. Based on this foundation, we discuss the formulation of optimization-based utility elicitation with stochastic dominance constraints for large-scale multi-attribute group decision-making and robust satisficing extensions to deal with uncertainty.
\end{itemize}

The remainder of this paper is organized as follows: Section \ref{sec-2} details the derivation of OPA and proves several fundamental properties. Section \ref{sec-3} introduces the unified framework, formulation, and validation metrics of GOPA. Section \ref{sec-4} demonstrates and validates GOPA through a case study on the improvisational emergency supplier selection during the 7.20 mega-rainstorm disaster in Zhengzhou, China. Section \ref{sec-5} discusses the formulation considering the elicitation errors, as well as the advantages and insights of GOPA. Section \ref{sec-6} presents the conclusions and outlines future research directions.

\section{Ordinal Priority Approach}
\label{sec-2}

\subsection{Derivation of Ordinal Priority Approach}
\label{subsec-2-1}

Consider a classical MADM problem that decision-maker needs to determine the optimal alternative from $K$ alternatives, $\mathcal{M}=\{1,2,\dots, K\}$, based on their performance on $J$ attributes $\mathcal{N}=\{1,2,\dots, J\}$, as assessed by $I$ experts $\mathcal{Q}=\{1,2,\dots, I\}$. Suppose that the independence axiom holds for both attributes and alternatives. Initially, decision maker assigns the ranking ${{t}_{k}}\in [I]$ to each expert $k\in \mathcal{Q}$ based on factors such as work experience, educational background, and organizational structure. Subsequently, each expert $k\in \mathcal{Q}$ independently provides the ranking ${{s}_{jk}}\in [J]$ of each attribute $j\in \mathcal{N}$ and the ranking ${{r}_{ijk}}\in [K]$ of each alternative $i\in \mathcal{M}$ on attributes $\forall j\in \mathcal{N}$. In OPA, each expert acts as an independent decision-maker, providing evaluations independently without needing group deliberation. Thus, all the rankings can reflect the individual preference information. Following customary conventions, the most crucial rank is designated as 1, the subsequent one as 2, and so on. Let $A_{ijk}^{({{r}_{ijk}})}$ represent the alternative $i$ with the ranking of ${{r}_{ijk}}$ on the attribute $j$ under the preferences of the expert $k$. 

To simplify the formulation in the following, we first define the following three sets:
\[
\begin{aligned}
& {{\mathcal{X}}^{1}} := \left\{(i,j,k,l)\in \mathcal{Q}\times \mathcal{N}\times \mathcal{M}\times \mathcal{M} : {{r}_{ijl}}={{r}_{ijk}}+1,{{r}_{ijk}}\in [K-1] \right\}, \\ 
& {{\mathcal{X}}^{2}} := \left\{(i,j,k)\in \mathcal{Q}\times \mathcal{N}\times \mathcal{M}:{{r}_{ijk}}=K \right\}, \\ 
& \mathcal{Y} := \left\{(i,j,k)\in \mathcal{Q}\times \mathcal{N}\times \mathcal{M} \right\}. \\ 
\end{aligned}
\]

For the same expert $i$ and attribute $j$, the alternatives with better rankings dominate those with worse rankings, as expressed in Equation (\ref{eq-01}).
\begin{equation}
A_{ijl}^{({{r}_{ijl}})} \preceq A_{ijk}^{({{r}_{ijk}})} \quad \forall (i,j,k,l)\in {{\mathcal{X}}^{1}}
\label{eq-01}
\end{equation}

Let ${{w}_{ijk}}$ denote the weight of alternative $i$ on attribute $j$ under the preferences of expert $k$. The normalized weight space is then defined as
\[
\mathcal{W} := \left\{{{w}_{ijk}} \in \mathbb{R}^{I \times J \times K}:\sum\limits_{i=1}^{I}{\sum\limits_{j=1}^{J}{\sum\limits_{k=1}^{K}{{{w}_{ijk}}}}}=1,{{w}_{ijk}}\ge 0,\forall i\in \mathcal{Q},j\in \mathcal{N},k\in \mathcal{M} \right\}.
\]

By the expected utility theory, the following statement holds:
\begin{equation}
A_{ijl}^{({{r}_{ijl}})}\preceq A_{ijk}^{({{r}_{ijk}})}\Leftrightarrow {{w}_{ijl}}\le {{w}_{ijk}}, \quad \forall (i,j,k,l)\in {{\mathcal{X}}^{1}}.
\label{eq-02}
\end{equation}

To incorporate the impact of expert preference on the weight disparity of each pair of alternatives with consecutive rankings, OPA multiplies both sides of the inequality in Equation (\ref{eq-02}) by the ranking parameters:
\begin{equation}
\begin{aligned}
\delta (\mathbf{w})=\left\{ \begin{array}{*{35}{l}}
   {{\delta }^{1}}({{w}_{ijk}})={{t}_{i}}{{s}_{ij}}{{r}_{ijk}}({{w}_{ijk}}-{{w}_{ijl}})\ge 0, \quad &&\forall (i,j,k,l)\in {{\mathcal{X}}^{1}},  \\
   {{\delta }^{2}}({{w}_{ijk}})={{t}_{i}}{{s}_{ij}}{{r}_{ijk}}({{w}_{ijk}})\ge 0, &&\forall (i,j,k)\in {{\mathcal{X}}^{2}}.  \\
\end{array} \right.
\end{aligned}
\label{eq-03}
\end{equation}

Equation (\ref{eq-03}) describes the marginal utility of weight increments derived from alternative ranking. In the normalized weight space, any point within the polyhedron satisfying the conditions of Equation (\ref{eq-03}) aligns with all experts’ preference rankings. Therefore, some studies offer a uniform distribution of marginal effects of all weight increments to accommodate every expert preference without adding extra information \citep{S22}. However, in the practical decision-making process, decision-makers are inclined to pursue decision-weight calculations that reflect the preferences of experts while demonstrating maximum differentiation \citep{W.S.C.S.C.G24}. Thus, in OPA, a multi-objective optimization model for decision-weight elicitation based on the ranking preference is given by
\begin{equation}
\max_{\mathbf{w}\in \mathcal{W}} \delta (\mathbf{w}).
\label{eq-04}
\end{equation}

By employing the max-min method and variable substitution, the multi-objective optimization model in Equation (\ref{eq-04}) can be transformed into linear programming, as depicted in Equation (\ref{eq-05}).
\begin{equation}
\begin{aligned}
\max_{\mathbf{w} \in \mathcal{W}, z} \text{ } & z \\ 
 \text{s.t. } & z\le {{\delta }^{1}}({{w}_{ijk}}) \quad && \forall (i,j,k,l)\in {{\mathcal{X}}^{1}} \\ 
 & z\le {{\delta }^{2}}({{w}_{ijk}}) && \forall (i,j,k)\in {{\mathcal{X}}^{2}}
\end{aligned}
\label{eq-05}
\end{equation}

\begin{proposition}
Given the rankings of experts $t_{i}, \forall i \in \mathcal{Q}$, attributes provided by experts $s_{ij}, \forall i \in \mathcal{Q}, \forall j \in \mathcal{N}$, and alternatives under attributes provided by experts $r_{ijk}, \forall i \in \mathcal{Q}, j \in \mathcal{N}, k \in \mathcal{M}$, the decision-weight elicitation model based on the ranking preference information can be formulated as Equation (\ref{eq-06}).
\begin{equation}
\begin{aligned}
\max_{\mathbf{w}, z} \text{ } & z \\ 
\mathrm{s.t.} \text{ } & z\le {{t}_{i}}{{s}_{ij}}{{r}_{ijk}}({{w}_{ijk}}-{{w}_{ijl}}) \quad &&\forall (i,j,k,l)\in {{\mathcal{X}}^{1}} \\ 
 &  z\le {{t}_{i}}{{s}_{ij}}{{r}_{ijk}}({{w}_{ijk}}) && \forall (i,j,k)\in {{\mathcal{X}}^{2}} \\ 
 &  \sum\limits_{i=1}^{I}{\sum\limits_{j=1}^{J}{\sum\limits_{k=1}^{K}{{{w}_{ijk}}}}}=1 \\ 
 &  {{w}_{ijk}}\ge 0 &&  \forall (i,j,k)\in \mathcal{Y}
\end{aligned}
\label{eq-06}
\end{equation}
\label{proposition-01}
\end{proposition}

The decision-weight elicitation model involves $I\times J\times K+1$ variables and $I\times J\times K+2$ constraints without counting the non-negativity constraints. In the case of individual decision-making, the decision-weight elicitation model can be formulated without the parameters and constraints related to multiple decision-makers.

After solving Equation (\ref{eq-06}), the weight of experts, attributes, and alternatives, denoted as $W_{i}^{\mathcal{Q}}$, $W_{j}^{\mathcal{N}}$, and $W_{k}^{\mathcal{M}}$, respectively, can be expressed by Equation (\ref{eq-07}).
\begin{equation}
\begin{aligned}
& W_{i}^{\mathcal{Q}}=\sum\limits_{j=1}^{J}{\sum\limits_{k=1}^{K}{{{w}_{ijk}}}} \quad && \forall i\in \mathcal{Q} \\ 
 & W_{j}^{\mathcal{N}}=\sum\limits_{i=1}^{I}{\sum\limits_{k=1}^{K}{{{w}_{ijk}}}} && \forall j\in \mathcal{N} \\ 
 & W_{k}^{\mathcal{M}}=\sum\limits_{i=1}^{I}{\sum\limits_{j=1}^{J}{{{w}_{ijk}}}} && \forall i\in \mathcal{M} \\
\end{aligned}
\label{eq-07}
\end{equation}

\subsection{Properties of Ordinal Priority Approach}
\label{sec-2-2}

In this section, we formally establish fundamental properties of OPA, including optimal solution efficiency, analytical solutions, and decomposability. These form the basis for analyzing the mathematical theory of OPA. We begin by evaluating the efficiency of the optimal solution derived from Proposition \ref{proposition-01} using the max-min method.

\begin{definition}
A solution $\mathbf{w}$ is said to dominate another solution $\mathbf{w}'$ if and only if there exists a specific element $(i,j,k)' \in \mathcal{Y}$ such that $\delta(\mathbf{w}) > \delta(\mathbf{w}')$, and $\delta(\mathbf{w}) \ge \delta(\mathbf{w}')$ for every $(i,j,k) \in \mathcal{Y} \backslash \{(i,j,k)' \}$.
\end{definition}

\begin{definition}
A feasible solution $\mathbf{w}$ is called efficient if it cannot be dominated by other solution $\mathbf{w}'$.
\end{definition}

Let us consider the two-stage mathematical programming, where Equation (\ref{eq-06}) from Proposition \ref{proposition-01} constitutes the first stage, with its optimal value denoted as $z^{*}$, followed by Equation (\ref{eq-08}) as the second stage.
\begin{equation}
\begin{aligned}
    z^{\Bar{*}} = \max\limits_{\mathbf{w}} \text{ } & \sum\limits_{i \in \mathcal{Q}} \sum\limits_{j \in \mathcal{N}} \sum\limits_{k \in \mathcal{M}} \delta(\mathbf{w}) \\
    \text{s.t. } & z^{*} \le \delta(\mathbf{w}) \quad &&\forall (i,j,k)\in {{\mathcal{X}}^{1}} \cup {{\mathcal{X}}^{2}} \\ 
    &  \sum\limits_{i=1}^{I}{\sum\limits_{j=1}^{J}{\sum\limits_{k=1}^{K}{{{w}_{ijk}}}}}=1 \\ 
    &  {{w}_{ijk}}\ge 0 &&  \forall (i,j,k)\in \mathcal{Y} 
    \label{eq-08}
\end{aligned}
\end{equation}

Let $OS(P)$ and $OS(Q)$ represent the optimal solution sets of the first-stage and second-stage problems, respectively, while $FS(P)$ and $FS(Q)$ denote the feasible solution sets of the first-stage and second-stage problems, respectively. Then, we have the following statement.

\begin{lemma}
    $z^{*} = \min \delta(\mathbf{w}) $ for $\mathbf{w} \in OS(Q)$.
    \label{lemma-01}
\end{lemma}

\begin{theorem}
    The solution obtained from the above two-stage mathematical programming is efficient.
    \label{theorem-01}
\end{theorem}

The above statement indicates that the optimal solution set derived from Proposition \ref{proposition-01} using the max-min method may not all be efficient. However, at least one element of the optimal solution set is an efficient point for the multi-objective optimization model delineated in Equation (\ref{eq-03}). Notably, the max-min solution is typically more manageable than the entire Pareto optimal set \citep{B.C24}.

In practical decision-making process, decision-makers may encounter situations where experts assign identical rankings to several alternatives or where missing values arise due to expert tendencies to exclude specific alternatives. Such scenarios lead to the structural alterations of OPA, particularly concerning constraints of weight disparities among alternatives with consecutive rankings. Thus, we initially map the alternative index $k$ in $w_{ijk}$ to the ranking index $r$ corresponding to their positions in ranking $r_{ijk}$ and define the set
\[
\mathcal{U} := \left\{ (i,j,r) : i\in [I],j\in [J],r\in [{{K}_{ij}}],{{K}_{ij}}=\underset{k}{\mathop{\max }}\,\{{{r}_{ijk}}\} \right\}.
\]

For expert $i$ and attribute $j$, there are ${{K}_{ij}}$ constraints in Equation (\ref{eq-06}), consisting of ${{K}_{ij}}-1$ constraints of $\delta^{1}({{w}_{ijr}})$ and 1 constraint of $\delta^{2}({{w}_{ijr}})$. Then, the cumulative sum of the last $r$ constraints in ascending order yields
\[
{{w}_{ijr}} =\frac{1}{t_{i}{s}_{ij}} \left(\sum\nolimits_{h=r}^{{{K}_{ij}}}{\frac{1}{h}} \right) z, \quad \forall (i,j,r)\in \mathcal{U}.
\]

Let ${{c}_{ijr}}$ denotes the frequency of ranking $r$ occurs on attribute $j$ given by expert $i$, and there exists $\sum\nolimits_{r=1}^{{{K}_{ij}}}{{{c}_{ijr}}}=K$ and $\sum\nolimits_{r=1}^{{{K}_{ij}}}{({{c}_{ijr}}-1)}=K-{{K}_{ij}}$. The weight assigned to the alternatives with the missing and duplicate rankings under the same expert and attribute should be 0 or the same, respectively. Thus, the OPA model in Proposition \ref{proposition-01} can be reformulated as Equation (\ref{eq-09}).
\begin{equation}
\begin{aligned}
\max_{\mathbf{w}, z} \text{ } & z \\
\text{s.t. } & \left( \sum\nolimits_{h=r}^{{{K}_{ij}}}{\frac{1}{h}} \right)z\le {{t}_{i}}{{s}_{ij}}{{w}_{ijr}} \quad && \forall (i,j,r)\in \mathcal{U} \\ 
& \sum\limits_{i=1}^{I}{\sum\limits_{j=1}^{J}{\sum\limits_{r=1}^{{{K}_{ij}}}{{{c}_{ijr}}{{w}_{ijr}}}}}=1 \\ 
& {{w}_{ijr}}\ge 0 \quad && \forall (i,j,r)\in \mathcal{U}  
\end{aligned}
\label{eq-09}
\end{equation}

The following theorem provides the analytical solution (closed-form solution) for OPA, enabling direct computation of optimal weights from this analytical solution rather than solving the LP problem.

\begin{theorem} 
The analytical solution of OPA is 
\begin{equation}
z^{*} = 1\Bigg/ \sum\limits_{i=1}^{I}{\sum\limits_{j=1}^{J}{\sum\limits_{r=1}^{{{K}_{ij}}}{\frac{{{c}_{ijr}}\sum\nolimits_{h=r}^{{{K}_{ij}}}{\frac{1}{h}}}{{{t}_{i}}{{s}_{ij}}}}}} 
\label{eq-10}
\end{equation} 
and 
\begin{equation}
 w^{*}_{ijr}= \left(\sum\nolimits_{h=r}^{{{K}_{ij}}}{\frac{1}{h}} \right)\Bigg/ \left({{t}_{i}}{{s}_{ij}}\sum\limits_{i=1}^{I}{\sum\limits_{j=1}^{J}{\sum\limits_{r=1}^{{{K}_{ij}}}{\frac{{{c}_{ijr}}\sum\nolimits_{h=r}^{{{K}_{ij}}}{\frac{1}{h}}}{{{t}_{i}}{{s}_{ij}}}}}} \right), \quad \forall (i,j,r)\in \mathcal{U}.
 \label{eq-11}
\end{equation}
\label{theorem-02}
\end{theorem}

After having the optimal solution, the weight of ${{w}_{ijr}}$ is assigned to the alternative ranked $r$ under attribute $j$ and expert $i$. By Theorem \ref{theorem-01}, when there are no missing or same rankings (i.e., ${{K}_{ij}}=K$ for every $i\in \mathcal{Q},j\in \mathcal{N}$), the analytical solution reduces to 
\begin{equation}
z^{*}=1\Bigg/\left(K \left(\sum\limits_{p=1}^{I}{\frac{1}{p}} \right) \left(\sum\limits_{q=1}^{J}{\frac{1}{q}} \right)\right)
\label{eq-12}
\end{equation}
and 
\begin{equation}
{{w}^{*}_{ijr}}= \left(\sum\limits_{h=r}^{K}{\frac{1}{h}}\right)\Bigg/\left({{t}_{i}}{{s}_{ij}}K \left(\sum\limits_{p=1}^{I}{\frac{1}{p}}\right) \left(\sum\limits_{q=1}^{J}{\frac{1}{q}} \right)\right),\quad \forall (i,j,r)\in \mathcal{U}.
\label{eq-13}
\end{equation}
In the following of this paper, we primarily utilize Equations (\ref{eq-12}) and (\ref{eq-13}) for property analysis, focusing on scenarios without missing or duplicate rankings.

\begin{example}
    Consider the case, including 3 experts, 5 attributes, and 10 alternatives without missing and duplicate rankings. Figure \ref{fig-01} illustrates the alternative weights with different rankings under experts and attributes and the aggregated weights of attributes and experts.
    \begin{figure}[h]
    \centering
    \includegraphics[width=\textwidth]{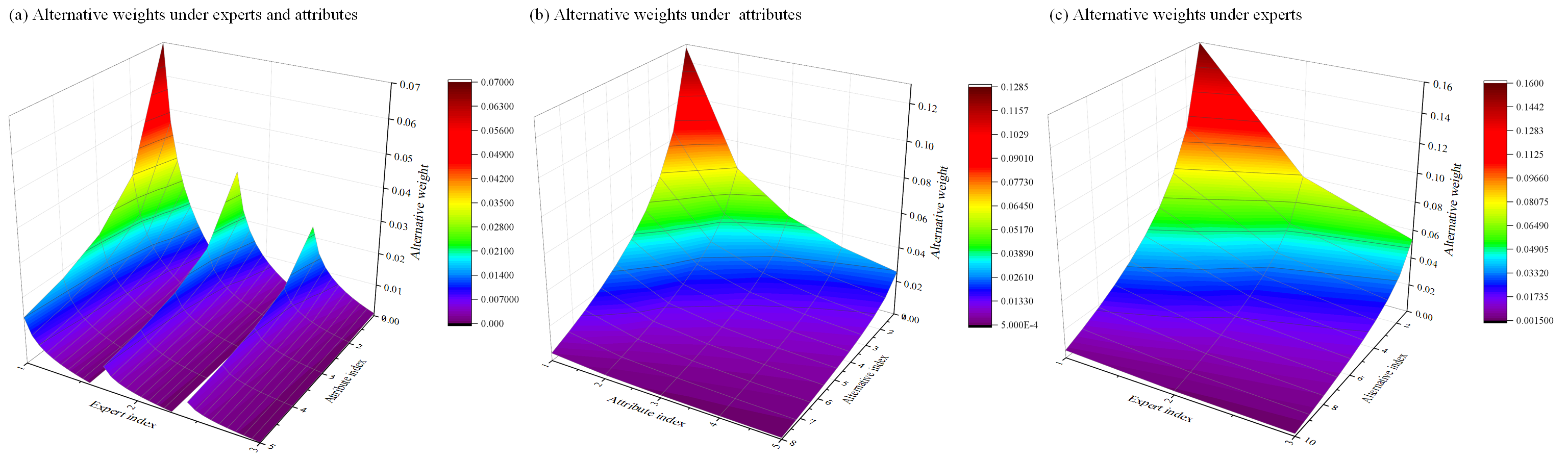}
    \caption{Weight results of the scenario with 3 experts, 5 attributes, and 10 alternatives}
    \label{fig-01}
    \end{figure}
    \label{example-01}
\end{example}

The results show that as the ranking increases, there is a diminishing marginal effect on the weight disparities of the alternatives with consecutive rankings. This characteristic is prevalent in rank-based surrogate weights, such as rank order centroid weights and rank reciprocal weights \citep{B.N23}. Therefore, we further explore the relationship between the analytical solution of OPA and rank-based surrogate weights.

\begin{corollary}
    Let $v_{r}^{ROC}$ represent the rank order centroid weight of the alternative ranked $r$, defined as $ v_{r}^{ROC} = \left(\sum_{h=r}^{K}\frac{1}{h}\right) \Big/ K$. Then, for every $(i,j,r)\in \mathcal{U}$, the optimal weight $w_{ijr}^{*}$ of the alternative ranked $r$ under expert $i$ and attribute $j$ in OPA can be expressed as
    \[
    w_{ijr}^{*} = \frac{v_{r}^{ROC}}{t_{i} s_{ij}\left(\sum_{p=1}^I\frac{1}{p}\right)\left(\sum_{q=1}^J\frac{1}{q}\right)}, \quad\forall(i,j,r)\in\mathcal{U},
    \]
    which implies that the rank order centroid weight is a special case of OPA for determining the weights of alternatives when the importance of experts and attributes is equal.
    \label{corollary-01}
\end{corollary}

\begin{corollary}
    Let $v_{s}^{RR}$ denote the rank reciprocal weight of the attribute ranked $s$ under expert $i$, which is given by $ v_{s}^{RR} = 1 \Big/ \left( s_{ij} \sum_{q=1}^J\frac{1}{q} \right)$. Then, for every $(i,s)\in \mathcal{U}$, the optimal weight $w_{is}^{*}$ of the attribute ranked $s$ under expert $i$ in OPA can be determined by
    \[
    w_{is}^{*} = \frac{v_{s}^{RR}}{t_{i}\left(\sum_{p=1}^{I} \frac{1}{p}\right)}, \quad \forall (i,s) \in \mathcal{U},
    \]
    which implies that the rank reciprocal weight is a special case of OPA for determining the weights of attributes when the importance of experts is equal.
    \label{corollary-02}
\end{corollary}

The following confirms the decomposability of the analytical solution in Theorem \ref{theorem-01}, a provides the foundation for extending OPA based on prospect theory. 

\begin{definition}
    A utility $u$ is called rank-based net utility if its value is solely determined by its ranking position $r$, i.e., there exists $f: \mathbb{R} \mapsto \mathbb{R}$ such that $u = f(r)$.
\end{definition}

\begin{proposition}
    For every $(i,r) \in \mathcal{U}$, the optimal weight $w_{ir}^{*}$ of the alternative ranked $r$ under expert $i$ in OPA can be decomposed into the product of the weight $w_{i}^{\mathcal{Q}}$ of expert $i$ and the rank-based net utility $u_{ir}$ of the alternative ranked $r$ under expert $i$, i.e., $w_{ir}^{*} = w_{i}^{\mathcal{Q}}u_{ir}, \forall (i,r) \in \mathcal{U}$.
    \label{proposition-02}
\end{proposition}

Proposition \ref{proposition-02} forms the foundation for extending OPA based on prospect theory. The rank-based net utility of the alternative provides a comparable unweighted outcome with a consistent scale across experts, which can serve as a reference point for evaluating gains and losses in prospect theory. Meanwhile, the weights of the experts can be regarded as the importance/probability of the corresponding prospects. And we leave the extension of OPA under prospect theory as a future direction to explore.

\section{Generalized Ordinal Priority Approach}
\label{sec-3}

In this section, we introduce the Generalized Ordinal Priority Approach (GOPA) from the perspective of the two-stage “estimate-then-optimize” contextual optimization. It integrates partial preference information with two mainstream utility structures, including the discrete rank-based surrogate weight and the continuous risk preference utility function.

\subsection{Unified Framework}
\label{sec-3-1}

 Recall that the left side of the inequality in the reformulated OPA model (Equation (\ref{eq-09})) can be interpreted as the product of the ROC weights and the number of alternatives. Thus, OPA can be regarded as eliciting utilities based on the ROC weights within a normalized decision space. This space covers experts, attributes, and alternatives in MADM, aiming to maximize weight disparities while adhering to ranking preferences. Naturally, we contemplate extending the ROC weights to a more general utility structure to elicit utilities within a normalized decision space. In this context, the utilities are treated as uncertainty parameters corresponding to the rankings of various alternatives under different experts and attributes. The partial preference information and prior utility structure provided by decision-makers can serve as auxiliary information for deriving alternative utilities \citep{S.C.D.F.F.V24}. To address uncertainty, this study utilizes a two-stage “estimate-then-optimize” contextual optimization to construct GOPA. Specifically, it first elicits the optimal utility distribution of alternatives across all attributes and experts based on partial preference information and utility structure, followed by optimizing to determine the optimal decision weights. This results in the unified GOPA framework, as expressed in Equation (\ref{eq-14}).
\begin{equation}
        \max\limits_{\mathbf{w} \in \mathcal{W} ,z} \left\{z : f(z,\mathbf{u}^{*}) \preceq g(\mathbf{w}), \mathbf{u}^{*} = \underset{\mathbf{u} \in \mathcal{P}}{\mathop{\arg \min }} \text{ } h(\mathbf{u}, \mathbf{v}) \right\}
    \label{eq-14}
\end{equation}

Where $f(z, \mathbf{u}^{*}_{ijr}) = {K}_{ij} u_{ijr}^{*} z \text{ and } g(\mathbf{w}_{ijr}) = {{t}_{i}}{{s}_{ij}}{{w}_{ijr}} \text{ for } \forall (i,j,r)\in \mathcal{U}$, $\mathbf{v}$ represents the target utility, $\mathcal{P}$ denotes the feasible region of $\mathbf{u}$ derived from partial preference information, and $h(\cdot, \cdot)$ denotes a general loss function.

This study employs cross-entropy utility minimization to elicit the optimal alternative utility distribution in the first-stage problem of OPA. The minimum cross-entropy utility theorem is proposed by \cite{A06} based on the utility density function concept, drawing an analogy between utility and probability. Entropy utility adheres to the axioms of Von Neumann and Morgenstern's expected utility theory while satisfying the fundamental independence requirement between utility and probability. Given the normative division between belief and preference, the utility value of a prospect remains unaffected by the probability of its realization \citep{J.N.W08}. In GOPA, the adoption of cross-entropy utility minimization is motivated by its adaptability to input data of varying statistical quality. This permits experts to furnish decision information that is sufficiently determined without strict adherence to statistical data quantity requirements, thereby aligning to alleviate the burden of expert decision-making. On the other hand, in the absence of prior structural information, cross-entropy utility minimization reduces to entropy utility maximization, facilitating the elicitation of utilities solely based on available partial preference information. It tackles the common challenge of accurately capturing the specific utility structures of decision-makers.

\begin{figure}[h]
\centering
\includegraphics[width=5.5in]{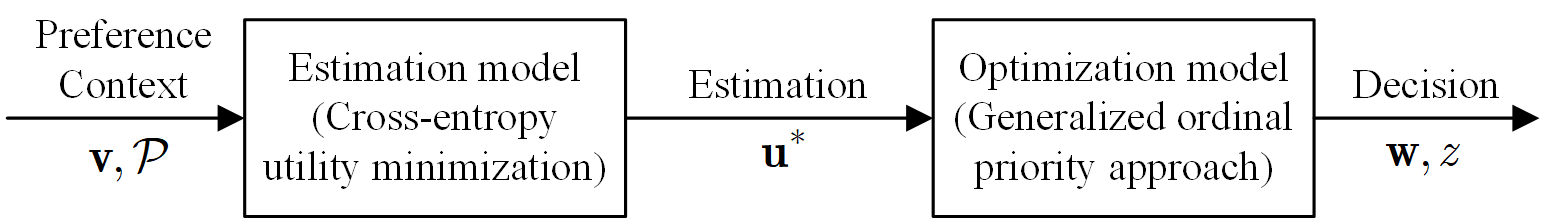}
\caption{GOPA decision pipeline under the partial preference information}
\label{fig-02}
\end{figure}

The GOPA decision pipeline is shown in Figure \ref{fig-02}. 

\subsection{Decision Information for Utility Elicitation}
\label{sec-3-2}
\subsubsection{Partial Preference Information}
\label{sec-3-2-1}

 The main types of partial preference information for utility elicitation include weak ordered relations, ratio scales, absolute differences, and lower bounds \citep{A15, A17}. In GOPA, the partial preference context for specific experts and attributes can be derived from the intersections of these sets of partial preference information:
\[
\begin{aligned}
    \mathcal{P}^{RS} & := \left\{ \mathbf{u} \in \mathbb{R}^{K} : u_{r}  = \alpha_{e_{1}} u_{r+1}, r \in [K], \forall e_{1} \in [E_{1}] \right\}, \\
    \mathcal{P}^{AD} & := \left\{ \mathbf{u} \in \mathbb{R}^{K} : u_{r} - u_{r+1} = \beta_{e_{2}}, r \in [K], \forall e_{2} \in [E_{2}] \right\}, \\
    \mathcal{P}^{LB} & := \left\{ \mathbf{u} \in \mathbb{R}^{K} : u_{r} \ge \gamma_{e_{3}} \ge 0, r \in [K], \forall e_{3} \in [E_{3}] \right\}, \\
    \mathcal{P}^{WO} & := \left\{ \mathbf{u} \in \mathbb{R}^{K} : u_{r} - u_{r+1}\ge 0, \forall r \in [K] \right\}, \\
    \mathcal{P}^{NN} & := \left\{ \mathbf{u} \in \mathbb{R}^{K} : u_{r} \ge 0, \forall r \in [K] \right\}, \\
    \mathcal{P}^{NZ} & := \left\{ \mathbf{u} \in \mathbb{R}^{K} : \sum\nolimits_{r=1}^{K} u_{r} = 1 \right\},
\end{aligned}
\]
where the utility of the alternative ranked $r$ is denoted as $u_{r}$; the set of absolute differences is denoted as $\mathcal{P}^{AD}$, which is typically used by the quasi-distance-based MADM methods such as TOPSIS, VIKOR, and EDAS; the set of ratio scales is denoted as $\mathcal{P}^{RS}$, which is commonly employed by pairwise comparison methods such as AHP, ANP, and BWM; the set of lower bounds is denoted as $\mathcal{P}^{LB}$; the set of weak ordered relations is denoted as $\mathcal{P}^{WO}$. $\mathcal{P}^{NN} \text{ and } \mathcal{P}^{NZ}$ represents the non-negative and normalized scaling conditions, respectively, which are typical conditions in expected utility theory. 

When $\beta_{e_1} = 0 \text{ and } \alpha_{e_{2}} = 1$, $\mathcal{P}^{AD} \text{ and } \mathcal{P}^{RS}$ reduce to $\mathcal{P}_{WO}$. This suggests that weak ordered relations (i.e., dominance ranking preference) can be foundational for integrating other reliable partial preference information (e.g., $\mathcal{P}^{AD}, \mathcal{P}^{RS},\text{ and } \mathcal{P}^{LB}$) from decision-makers to approximate actual preferences. Notably, experts need only to provide sufficiently determined preference information, implying that the mapping of partial preference information indexes to alternative indexes is not always surjective. This forms the basis for personalized preference analysis under different experts and attributes in GOPA. Consequently, for the specific expert $i$ and attribute $j$, the partial preference context can be expressed as $\mathcal{P}_{ij} =  \mathcal{P}^{NN} \cap \mathcal{P}^{NP} \cap \mathcal{P}^{WO} \cap \mathcal{P}^{AD}_{ij} \cap \mathcal{P}^{RS}_{ij} \cap \mathcal{P}^{LB}_{ij}$ across experts and attributes. When $\mathcal{P}^{AD}_{ij}, \mathcal{P}^{RS}_{ij},\text{ and } \mathcal{P}^{LB}_{ij}$ do not exists, preference context will reduce to the unbiased context.

\subsubsection{Global Utility Structure}
\label{sec-3-2-3}

The global utility structure in GOPA forms the prior structure information for eliciting the alternative utility. This study provides the global utility structure in discrete and continuous prospects. The utility structure in the discrete form includes the rank-based surrogate weights, and the continuous form involves the risk preference utility function. Notably, the absence of global utility information does not affect the implementation of GOPA.

As for the discrete form of global utility structure, rank-based surrogate weighting is the primary method based on dominance preference for eliciting attribute utility \citep{L.L.L20}. This approach can facilitate the objective reporting of subjective observations and experiences, mitigating response errors and judgmental bias, particularly when precise numerical information is lacking. The following outlines the prevalent methods of rank-based surrogate weights \citep{B.N23}. Notably, the need to elicit utility via rank-based surrogate weights aligns with that of $\mathcal{P}^{WO}$, avoiding additional decision-making burdens.

The rank sum (RS) weight is derived by standardizing individual ranks by the sum of ranks, with the weight being proportional to its ranking, which is given by
\[
    v_r^{RS}=\frac{K+1-r}{\sum_{t=1}^Kt}=\frac{2(K+1-r)}{K(K+1)},\quad\forall r\in[K].
\]

The rank exponent fixed (REF) weight is derived by raising the ranks obtained from the RS weight to the power of $z$:
\[
    v_r^{REF}=\frac{(K+1-r)^z}{\sum_{t=1}^Kt^z},\quad\forall r\in[K],
\]
where $z = 1.17$ is a common value in numerous studies \citep{B.N23}. 

The rank reciprocal (RR) weight defines the weight of each attribute as the reciprocal of its rank, then normalizes the ranks by dividing by the sum of all reciprocals:
\[
    v_r^{REF}=\frac{1}{r\sum_{j=1}^K\frac{1}{j}},\quad\forall r\in[K].
\]

The sum reciprocal (SR) weight is derived by combining the RS and RR weights. It averages the weights generated by these two methods to produce a distribution with a smaller range than either individual weight:
\[
    v_r^{REF}=\frac{\left(\frac{K+1-r}{K}+\frac{1}{r}\right)}{\sum_{j=1}^K\left(\frac{K+1-j}{K}+\frac{1}{j}\right)},\quad\forall r\in[K].
\]

The rank order centroid (ROC) weight identifies a set of weights representing all possible acceptable weight combinations:
\[
    v_r^{REF}=\frac{1}{K}\sum_{t=r}^K\frac{1}{t},\quad\forall r\in[K].
\]
This concept is based on the reasonable assumption that actual weights may be distributed anywhere within the feasible polytope, with its centroid coordinates being considered the optimal surrogates.

\begin{figure}[h]
\centering
\includegraphics[width=2.25in]{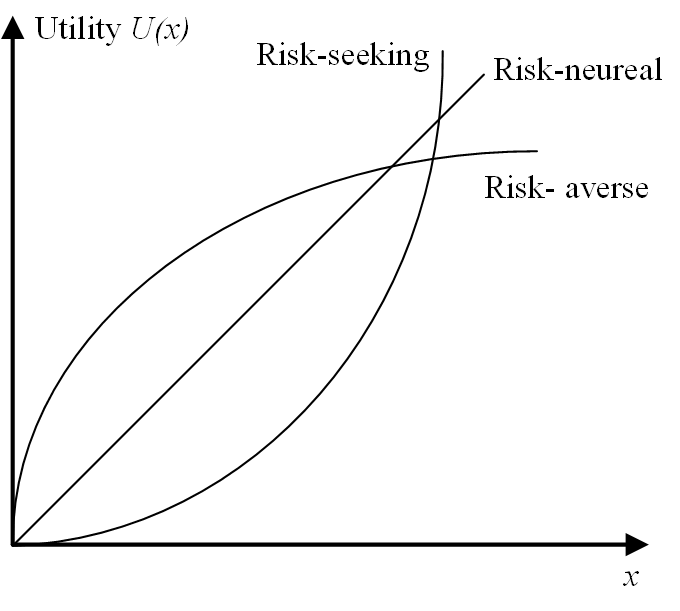}
\caption{Utility function and risk behaviour}
\label{fig-03}
\end{figure}

Regarding the continuous global utility structure in GOPA, risk preference utility functions are chosen. Abundant evidence demonstrates that risk preference significantly impacts practical decision-making processes, particularly in finance, economics, and decision science \citep{L.K.A23, C.S24}. It typically molds decision-makers’ cognitive processing of information, influencing their evaluation of alternatives and ultimately guiding their inclination towards optimal solutions \citep{P.D24}. Decision-makers’ risk preferences generally manifest as risk-seeking, risk-neutral, and risk-averse tendencies, which can be elucidated through utility functions. The utility functions representing risk-seeking, risk-neutral, and risk-averse behaviors exhibit concave, linear, and convex characteristics, respectively, as illustrated in Figure \ref{fig-03}.

Among the different risk preference utility functions, hyperbolic absolute risk aversion (HARA) has gained widespread application due to its generality \citep{X.W.Z24}. The HARA utility function is given by
\[
    V(x)=\frac1{1-\gamma}\bigg(\gamma\bigg(\beta+\frac\alpha\gamma x\bigg)^{1-\gamma}-1\bigg),
\]
with a utility density
\[
    v(x)=\alpha\left(\beta+\frac\alpha\gamma x\right)^{-\gamma},
\]
where $\alpha, \beta, \text{ and } \gamma$ are given constants. The HARA utility function can be generalized to common utility function forms:
\begin{itemize}
    \item[$-$] when $\gamma = 0$, it reduces to a risk-neutral utility function;
    \item[$-$] when $\beta = 0$ and $\gamma>0$, it reduces to a constant relative risk-averse (CRRA) utility function;
    \item[$-$]  when $\gamma \to \pm\infty$, it reduces to a constant absolute risk aversion (CARA) utility function.
\end{itemize}

Besides the above risk preference utility functions with single feature, the prospect theory proposed by \citet{T.K74} represents a significant approach to bridging the gap between normative theories and empirical behaviors. It suggests that decision-making is influenced by four factors: (1) outcome evaluation relative to a reference point; (2) losses have a more significant impact on decision-makers than gains; (3) increasing gains or losses weakens people's perceptions of them; (4) individuals exhibit probability distortion, tending to overweight small probabilities and underestimate large ones. These observations indicate that individuals exhibit risk aversion toward gains and risk-seeking toward losses. Thus, utility functions based on prospect theory exhibit an S-shaped characteristic. Consistent with normative research practices, the possibility of any probability evaluation bias are disregarded in this study, focusing instead on how to incorporate information indicating specific shapes of utility functions \citep{A.D15}. Therefore, this study adopts a logistic function as the utility function with the S-shaped characteristic, given by
\[
    V(x)=\frac1{1+e^{-x}}.
\]
We need further perform translation and stretching transformations to make it symmetric about point $\left((1+K_{ij})/2, 0.5\right)$.

Compared to rank-based surrogate weights, the risk preference utility function can complement the utility structure of GOPA in continuous prospects, which reveals the risk preference of the decision-makers.

\begin{remark}
    By the baseline assumption, the utility function must be continuous and twice-differentiable \citep{G.L17}. Despite not satisfying twice-differentiability, the rank-based surrogate weights after continuous transformation can be regarded to exhibit a risk-seeking tendency. Since, both sub-gradient and second-order sub-gradient of rank-based surrogate weights are greater than 0 (i.e., $\nabla v(r) > 0$ and $\nabla^{2} v(r) > 0$), aligning with the characteristics of risk-seeking utility function (i.e., $V^{'}(r) > 0$ and $V^{''}(r) > 0$).
    \label{remark-01}
\end{remark}

\subsection{Formulations}
\label{sec-3-3}

\subsubsection{Utility Elicitation in Discrete Prospect}
\label{sec-3-3-1}

In discrete prospects of GOPA, the utility structure of alternatives is derived from discrete rank-based surrogate weights as the target utility, incorporating partial preference context as constraints. Based on the cross-entropy utility minimization, Equation (\ref{eq-15}) expresses the utility elicitation in discrete prospects of the first-stage problem of GOPA.
\begin{subequations}\label{eq-15}
\begin{align}
    U_{ij}^{*}(r)=\underset{{\mathbf{U}_{ij}}(\cdot)}{\mathop{\arg \min }}\,\text{ } & -\sum\limits_{r=1}^{{{K}_{ij}}}{{{U}_{ij}}(r)\ln \frac{{{U}_{ij}}(r)}{V(r)}} \label{eq-15-1}\\ 
    \text{s.t. } & U_{ij}(r) - {\alpha_{ije_{1}}}(r)U_{ij}(r+1) = 0 \quad && \forall e_{1} \in [E_{1}] \label{eq-15-2}\\ 
    & {{U}_{ij}}(r)-{{U}_{ij}}(r+1)={{\beta }_{ije_{2}}}(r) \quad && \forall e_{2} \in [E_{2}] \label{eq-15-3}\\ 
    & {{U}_{ij}}(r)\ge {{\gamma }_{ije_{3}}} \quad &&\forall e_{3} \in [E_{3}] \label{eq-15-4}\\ 
    & {{U}_{ij}}(r)-{{U}_{ij}}(r+1) \ge 0 \quad && \forall r\in [{{K}_{ij}}] \label{eq-15-5}\\
    & \sum\limits_{r=1}^{{{K}_{ij}}}{{{U}_{ij}}(r)=1} \label{eq-15-6} \\
    & {{U}_{ij}}(r) \ge 0 \quad && \forall r\in [{{K}_{ij}}] \label{eq-15-7}
\end{align}
\end{subequations}

Where Equations (\ref{eq-15-2})-(\ref{eq-15-7}) indicate $\mathcal{P}^{RS}, \mathcal{P}^{AS}, \mathcal{P}^{LB}, \mathcal{P}^{WO}, \mathcal{P}^{NZ},\text{ and } \mathcal{P}^{NN}$, respectively. 

\begin{example}
    Consider 7 alternatives for ranking. Assuming that the lower bound of utilities is 0.03, the ratio scale difference of utilities between the alternative with the second and third rankings is 15$\%$, and the absolute preference difference of utilities between the alternative with the fourth and fifth rankings is 0.065. Figure \ref{fig-04} illustrates the comparison of rank-based surrogate weights with and without partial preference information derived from the cross-entropy utility minimization.
    \begin{figure}[h]
    \centering
    \includegraphics[width=4.5in]{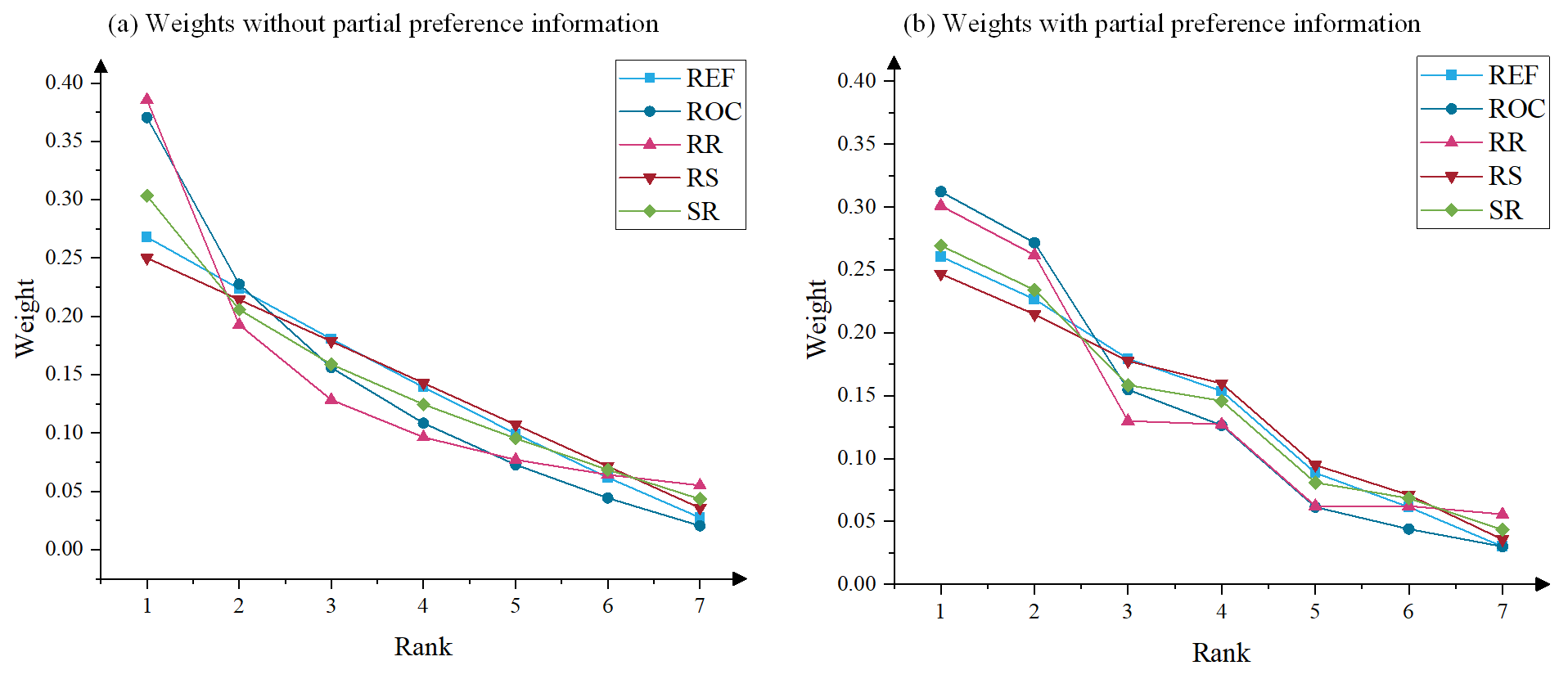}
    \caption{Comparison of rank-based surrogate weights with partial preference information}
    \label{fig-04}
    \end{figure}
\end{example}

In the absence of a pre-defined utility structure of rank-ordered surrogate weights, the following theorem shows that cross-entropy utility minimization will reduce to entropy utility maximization with partial preference context.

\begin{theorem}
    Let $U_{ij}^{EU}(r)$ represent the optimal utility of the alternative ranked $r$ under expert $i$ and attribute $j$ in the entropy utility maximization problem with partial preference context $\mathcal{P}_{ij}$ and let $\Bar{V}$ denote the uniformly weighted utility. Then, maximum entropy utility with partial preference context of the first-stage problem of GOPA in discrete prospects is a special case of minimum cross-entropy utility when the target utility form is uniformly weighted, expressed as: 
    \[
    U_{ij}^{EU}(r) = \underset{{\mathbf{U}_{ij}} \in \mathcal{P}_{ij}}{\mathop{\arg \max }} - \sum\limits_{r=1}^{K_{ij}}U_{ij}(r) \ln U_{ij}(r)  = \underset{{\mathbf{U}_{ij}} \in \mathcal{P}_{ij}}{\mathop{\arg \min }}  \sum\limits_{r=1}^{{{K}_{ij}}}{{{U}_{ij}}(r)\ln \frac{{{U}_{ij}}(r)}{\Bar{V}(r)}} .
    \]
    \label{theorem-03}
\end{theorem}

\subsubsection{Utility Elicitation in Continuous Prospect}
\label{sec-3-3-2}

 The optimal utility distribution for GOPA in continuous prospects is derived from the risk preference utility function as the target utility density function while incorporating partial preference context into the constraints. Equation (\ref{eq-16}) presents the elicitation formula for the utility density function under continuous prospects of the first-stage problem of GOPA.

\begin{subequations}\label{eq-16}
\begin{align}
     u^{*}_{ij}(x) = \underset{{{u}_{ij}}(x)}{\mathop{\arg \min }} \text{ } & \int_{0}^{K_{ij}} u_{ij}(x) \ln \left(\frac{u_{ij}(x)}{v(x)}\right) \mathrm{d} x \label{eq-16-1} \\
     \text{s.t. } & \int_{0}^{r} u_{ij}(x) \mathrm{d} x - \alpha_{ij{e}_{1}}\int_{0}^{r-1} u_{ij}(x) \mathrm{d} x = 0 && \quad \forall e_{1} \in [E_{1}] \label{eq-16-2} \\
     & \int_{0}^{r} u_{ij}(x) \mathrm{d} x - \int_{0}^{r-1} u_{ij}(x) \mathrm{d} x = \beta_{ij{e}_{2}} && \quad \forall e_{2} \in [E_{2}] \label{eq-16-3} \\
     & \int_{0}^{r} u_{ij}(x) \mathrm{d} x = \gamma_{ije_{3}} && \quad \forall e_{3} \in [E_{3}] \label{eq-16-4} \\
     & \int_{0}^{K_{ij}} u_{ij}(x) \mathrm{d} x = 1 \label{eq-16-5} \\
     & u_{ij}(x) \ge 0 \label{eq-16-6}
\end{align}
\end{subequations}

Where Equations (\ref{eq-16-2})-(\ref{eq-16-6}) represent $\mathcal{P}^{RS}, \mathcal{P}^{AS}, \mathcal{P}^{LB}, \mathcal{P}^{NZ},\text{ and } \mathcal{P}^{NN}$, respectively. Notably, $\mathcal{P}^{WO}$ is not contained since it can be guaranteed by the non-negativity of utility density function. 

\begin{example}
    Consider 7 alternatives for ranking. Assuming that the lower bound utility of the first-ranked alternative is 0.32, the ratio scale difference in utility between the second and third-ranked alternatives is 15$\%$, and the absolute utility difference between the fourth and fifth-ranked alternatives is 0.065. For the HARA utility function, $\alpha = 2, \beta = 1, \text{ and } \gamma = 1.5$, and for the CRRA utility function, $\alpha = 1 \text{ and } \gamma = 0.5$. Figure \ref{fig-05} presents a comparison of results among various risk preference utility functions. Notably, utility value increases monotonically with the ranking. 
    \begin{figure}[h]
    \centering
    \includegraphics[width=4.5in]{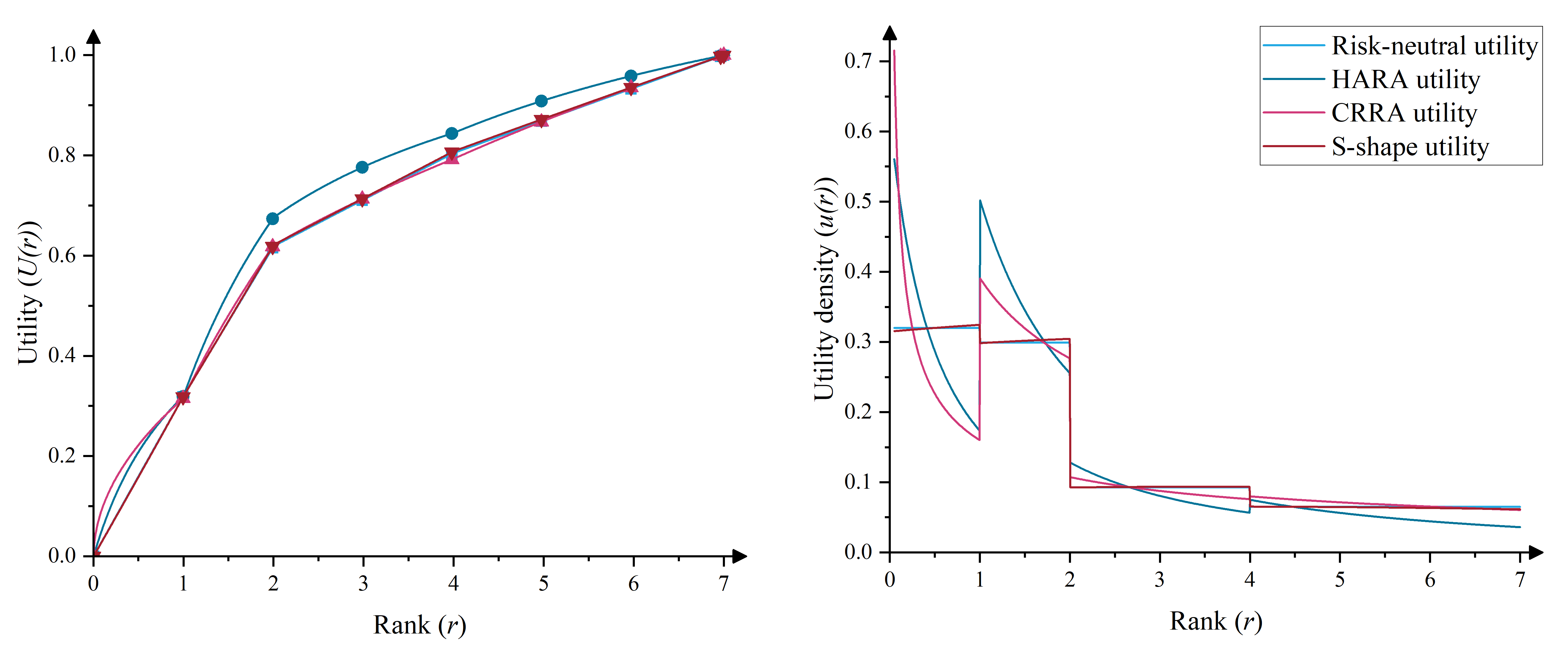}
    \caption{Comparison of utility function with partial preference information}
    \label{fig-05}
    \end{figure}
    \label{example-03}
\end{example}

The results reveal that the utility density function elicited by cross-entropy utility minimization exhibits a piecewise form, with breakpoints corresponding to the rankings based on the provided partial preference context. Integrating the elicited utility density functions yields utility functions with similar trends. Notably, despite the target utility density function's varying concavity and convexity features, the utility function derived from partial preference information does not display corresponding characteristics. The elicited utility density functions show that while utility densities exhibit concavity and convexity within different segments, they remain centered around a risk-neutral utility density function.  Thus, it can be inferred that in minimizing cross-entropy utility, partial preference information primarily shapes the utility function, while the target utility mainly adjusts within segments. In fact, this is a favorable phenomenon. As many studies have pointed out, determining decision-makers' risk preferences is challenging and prone to elicitation errors \citep{A.D15, G.X.Z23}. This phenomenon may somewhat diminish the influence of the global utility function, which benefits the cases when the global utility functions are absent. Based on the above findings, we proceed to analyze the risk preference relations between the elicited utility and the target utility density function of GOPA in continuous prospects with partial preference context.

\begin{proposition}
    Let $r_{e}, \forall e\in[E]$ represent the rankings corresponding to the partial preference context $\mathcal{P}_{ij}$ and let $v(x)$ denote the target utility density function. Then, the optimal solution $u_{ij}^{*}(x)$ to minimum cross-entropy utility with partial preference context of the first-stage problem of GOPA in continuous prospect is a piecewise function with breakpoints at $r_{e}, \forall e\in [E]$, which is given by
    \[
    u_{ij}^{*}(x) = v(x)e^{-1-\lambda_0-\sum_{e\in[E]}\lambda_e\zeta_e(x)},
    \]
    where $\zeta_e(x)$ is a step function with jumps at $r_{e}, \forall e\in [E]$.
    \label{proposition-03}
\end{proposition}

\begin{theorem}
    By Arrow-Pratt's definition of risk preference, the optimal solution $u_{ij}^{*}(x)$ to minimum cross-entropy utility with partial preference context of the first-stage problem of GOPA in continuous prospect has the identical risk preference $\eta_{ij}(x)$ with the target utility density function in each interval of $(0, r_1] \cup \left(\bigcup_{e\in[E]\backslash\{1\}}(r_{e-1},r_{e}]\right) \cup (r_{E},K_{ij}]$. Specifically, the risk preference within each interval is given by
    \[
    \eta_{ij}(x)=-\frac{\mathrm{d}}{\mathrm{d}x}\ln(u_{ij}^*(x))=-\frac{\mathrm{d}}{\mathrm{d}x}\ln(v(x)).
    \]
    \label{theorem-04}
\end{theorem}

\begin{corollary}
    Let $u_{ij}^{EU}(x)$ represent the optimal utility density function of expert $i$ and attribute $j$ in the entropy utility maximization problem with partial preference context $\mathcal{P}_{ij}$ and let $\Bar{v}(x)$ denote the uniform utility density function. Then, maximum entropy utility of GOPA in continuous prospects with partial preference context is a special case of minimum cross entropy utility when the target utility density function is uniform, expressed as:
    \[
    u_{ij}^{EU}(x) = \underset{{{u}_{ij}}(x) \in \mathcal{P}_{ij}}{\mathop{\arg \max }} - \int_{0}^{K_{ij}} u_{ij}(x) \ln \left(u_{ij}(x)\right) \mathrm{d} x = \underset{{{u}_{ij}}(x) \in \mathcal{P}_{ij}}{\mathop{\arg \min }} \int_{0}^{K_{ij}} u_{ij}(x) \ln \left(\frac{u_{ij}(x)}{\Bar{v}(x)}\right) \mathrm{d} x.
    \]
    \label{corollary-03}
\end{corollary}

\begin{remark}
    An alternative way to elicit utility for GOPA in continuous prospects with partial preference context is maximum entropy utility, which incorporates target utility density function into constraints. For example, let us consider the HARA utility function on a bounded domain: $V(x) = \frac1{1-\gamma}(\gamma(\beta+\frac\alpha\gamma x)^{1-\gamma}-1), x \in [0, K_{ij}]$. It has the utility density function, which can be rewritten into $v(x) = \alpha(\beta+\frac\alpha\gamma x)^{-\gamma}=e^{\ln(\alpha)-\gamma\ln(\beta+(\alpha/\gamma)x)}$, which is in accordance with the form of optimal solution shown in Equation (\ref{eq-sub43-03}). However, this formulation presents challenges in calculating the value of the Euler-Lagrange operators while optimizing the functional for entropy utility maximization. In contrast, the minimum cross-entropy utility utilized in this study exploits the properties of the indicator function, resulting in a utility density function in the form of a step function, eliminating the requirements to solve the Euler-Lagrange operators. Nonetheless, integrating the target utility density function into the constraints still holds promising prospects for depicting local risk preferences within a certain ranking segment, which we leave as a future direction to explore.
\end{remark}

The above discussion suggests that the risk preference of optimal solutions for minimizing cross-entropy utility can be controlled by selecting appropriate target utility functions. In cases where the target utility structure information is lacking, maximum entropy utility based on partial preference information can be employed to elicit utility, characterized by risk-neutral preferences.

\subsubsection{Decision Weight Optimization}
\label{sec-3-3-3}

After eliciting the utility of alternatives with specific rankings under different experts and attributes in discrete and continuous prospects (i.e., Equations (\ref{eq-15}) and (\ref{eq-16})), these are then integrated into the decision weight optimization problem to determine the optimal decision weights. Equation (\ref{eq-17}) expresses the decision weight optimization problem in the second-stage of GOPA.
\begin{equation}
    \begin{aligned}
        \max\limits_{\mathbf{w},z} \text{ } & z \\ 
     \text{s.t. } & {U^{*}_{ij}}(r)z \leq \frac{{{t}_{i}}{{s}_{ij}}{{w}_{ijr}}}{{K}_{ij}} \quad && \forall (i,j,r)\in \mathcal{U}^1 \\
     & \frac{\int_{0}^{r} u^{*}_{ij}(K_{ij} - x) \mathrm{d}x}{\sum\nolimits_{r=1}^{K_{ij}} \int_{0}^{r} u^{*}_{ij}(K_{ij} - x)\mathrm{d}x } z \leq \frac{{{t}_{i}}{{s}_{ij}}{{w}_{ijr}}}{{K}_{ij}} \quad && \forall (i,j,r)\in \mathcal{U}^2 \\ 
     & \sum\limits_{i=1}^{I}{\sum\limits_{j=1}^{J}{\sum\limits_{k=1}^{{{K}_{ij}}}{{{c}_{ijr}}{{w}_{ijr}}}}}=1 \\ 
     & {{w}_{ijr}}\ge 0 \quad && \forall (i,j,r)\in \mathcal{U} 
    \end{aligned}
    \label{eq-17}
\end{equation}
Where $\mathcal{U}^1$ and $\mathcal{U}^2$ represent the sets of indices for utility elicitation based on discrete and continuous prospects, respectively, and $\mathcal{U}^1 \cup \mathcal{U}^2 = \mathcal{U}$. Normalization of continuous prospect utility ensures alignment with the magnitude and direction of discrete prospect utility. Notably, this implies that GOPA can be extended to hybrid data sets with varying dimensions through suitable standardization techniques, which is a future direction to explore. After solving Equation (\ref{eq-17}), optimal weights $w^{*}_{ijr}$ are mapped to $w_{ijk}$ based on alternative rankings to assign weights to specifically ranked alternatives across experts and attributes. Finally, Equation (\ref{eq-07}) is utilized to calculate the weight of experts, attributes, and alternatives.

\begin{theorem}
    The analytical solution of the second-stage problem of GOPA in hybrid prospects is given by
    \begin{equation}
        z^{*} = 1 \Bigg/ \left(\sum_{i=1}^I\sum_{j=1}^J\sum_{r=1}^{K_{ij}}\frac{c_{ijr} K_{ij}U^{*}_{ijr}}{t_is_{ij}} \right),
        \label{eq-18}
    \end{equation}
    and
    \begin{equation}
        w^{*}_{ijr}= \left(K_{ij}U^{*}_{ijr}\right)\Bigg / \left(t_is_{ij}\sum_{i=1}^I\sum_{j=1}^J\sum_{r=1}^{K_{ij}}\frac{c_{ijr} K_{ij}U^{*}_{ijr}}{t_is_{ij}} \right),\quad\forall(i,j,r)\in\mathcal{U},
        \label{eq-19}
    \end{equation}
    where
    \[
    U_{ijr}^{*}= \begin{cases} U_{ij}^{*}(r), &\forall(i,j,r) \in \mathcal{U}^1, \\
    \left( \int_{0}^{r}u_{ij}^{*}(K_{ij}-x)dx \right) \Big/ \left( \sum_{r=1}^{K_{ij}}\int_{0}^{r}u_{ij}^{*}(K_{ij}-x)dx \right),  &\forall(i,j,r) \in \mathcal{U}^2.
    \end{cases}
    \]
    \label{theorem-05}
\end{theorem}

\begin{corollary}
    Let $\Gamma$ denote the set of all possible global utility structures, which can be either discrete or continuous. In the case without missing and duplicate rankings, the optimal objective value $z^{*}$ of the second-stage problem of GOPA is invariant for $\forall v \in \Gamma$.
    \label{corollary-04}
\end{corollary}

It follows from Corollary \ref{corollary-04} that the results of GOPA are robust. This stems from the fact that the assigned weight is determined by $w^{*}_{ijr} = (K_{ij}U^{*}_{ijr}z) / (t_{i}s_{ij})$, where the optimal value $z$ remains constant across various global utility structures. Hence, $w^{*}_{ijr}$ exhibits a linear relationship with $U^{*}_{ijr}$, which is constrained by the partial preference information and normalized scaling. On the other hand, the analytical solution provides a general lower-bound reference for other utility (weight) elicitation methods, which is not limited to GOPA. This forms the foundation for translating the subjective preference-based utility elicitation problem to stochastic dominance constraint problems, which we will discuss further in Section \ref{sec-5-1}.

\begin{remark}
    The unknown utility benchmark (lower bound) is the primal distinction of the subjective preference-based utility elicitation in this study compared to other utility elicitation problems, such as portfolio decisions. Consider the following optimization with stochastic dominance constraints: $\mathrm{max} \left\{ f(X)| Y \preceq_{(2)} X, X \in C  \right\}$ \citep{P.D24}. In portfolio decision problems, $Y$ can be viewed as a benchmark return, such as the performance of an existing portfolio or a market index. This ensures that no risk-neutral decision maker will prefer $Y$ over the optimal solution $X^{*}$ to the given problem \citep{D.R03}. However, such a benchmark is inaccessible for utility elicitation problems based on subjective preferences. Nevertheless, according to Corollary \ref{corollary-04}, the analytical solution in Theorem \ref{theorem-05} could provide an analogous benchmark reference for the general utility elicitation problem based on subjective preferences.
    \label{remark-03}
\end{remark}

\begin{corollary}
     In the case without missing and duplicate rankings, the weights of experts $w_{i}^{\mathcal{Q}}, \forall i \in \mathcal{Q}$ and attributes $w_{j}^{\mathcal{N}}, \forall j \in \mathcal{N}$ are independent of the experts' risk preferences toward the attributes.
     \label{corollary-05}
\end{corollary}

Corollary \ref{corollary-05} and Theorem \ref{theorem-05} demonstrate that in GOPA, experts' risk preferences towards distinct attributes only affect the weights assigned to alternatives, not the weights of experts and attributes, showcasing the property of risk preference independence of GOPA.

\subsection{Model Validation for Group Decision-Making}
\label{sec-3-4}

In this section, we introduce several statistical indicators to test the stability and reliability of GOPA group decision-making results, including percentage standard deviation, correlation coefficient, and confidence level measurement. Notably, the above metrics are optional and aid decision makers in understanding the state of group decision outcomes.

\subsubsection{Percentage Standard Deviation}
\label{sec-3-4-1}

Recall that the results of Example \ref{example-01} reveals a diminishing marginal effect in weight assignments across different rankings, leading to disparities between measurement scales. To mitigate this influence, dimensionless percentage standard deviation (PSD) are proposed to assess the dispersion of weights derived from expert preferences in GOPA \citep{M.J23}.

\begin{definition}
    The percentage standard deviation for attributes $\eta_{j}$ and alternatives $\eta_{k}$ in GOPA are defined by
    \begin{equation}
        \eta_{j} = \frac{1}{w^{\mathcal{N}}_j} \sqrt{\frac{\sum_{i=1}^I\left(w^{\mathcal{N}}_j \big/ I -w_{ij} \right)^2}{I-1}}, \quad \forall j \in \mathcal{N},
        \label{eq-20}
    \end{equation}
    and 
    \begin{equation}
        \eta_{k} = \frac{1}{w^{\mathcal{M}}_k} \sqrt{\frac{\sum_{i=1}^I\left(w^{\mathcal{M}}_k \big/ I -w_{ik} \right)^2}{I-1}}, \quad \forall k \in \mathcal{M},
        \label{eq-21}
    \end{equation}
    where $w^{\mathcal{N}}_j$ and $w^{\mathcal{M}}_k$ denote the weight of attribute $j$ and alternative $k$, respectively; $w_{ij}$ and $w_{ik}$ denote the weight of attribute $j$ and alternative $k$ derived from expert $i$, respectively.
    \label{definition-04}
\end{definition}

A higher PSD indicates a more dispersed distribution of weights, reflecting greater divergence in expert opinions. The weight outcomes are generally considered concentrated when PSD is below 0.2.

\begin{corollary}
    The percentage standard deviation for attributes $\eta_{j}, \forall j \in \mathcal{N}$ is independent of the experts' risk preference towards attributes.
    \label{corollary-06}
\end{corollary}

Notably, PSD for alternative weights is correlated with the risk preferences of experts towards attributes. However, the direction of change in PSD remains uncertain when risk preferences shift. For instance, if an expert's risk preference changes from risk-seeking to risk-averse, the weight allocation might appear more balanced intuitively, yet its relationship with the mean value remains indeterminate.

\subsubsection{Correlation Testing}
\label{sec-3-4-2}

Kendall's W correlation analysis, a traditional non-parametric method, is utilized to assess the consistency and significance of experts' decision outcomes in GOPA \citep{K.S23}. In GOPA, Kendall's W correlation analysis encompasses both alternatives and attributes with corresponding local and global confidence level measures. For alternatives, the decision outcomes of experts for each attribute are analyzed. Given attribute $j \in \mathcal{N}$, the weights $w_{ijk}$ assigned to alternative $k$ by the same expert $i$ are sorted in ascending order, denoted as $R_{ijk}$. When there are no duplicate rankings of alternatives, Kendall's W correlation coefficient of the alternatives under attribute $j$ is given by 
\begin{equation}
    \rho_{j}^{\mathcal{M}} = \frac{12 \sum_{k=1}^{K} (\sum_{i=1}^{I} R_{ijk} - \frac{1}{K} \sum_{i=1}^{I}\sum_{k=1}^{K} R_{ijk})^{2}}{I^{2}(K^{3}-K)}, \quad \forall j \in \mathcal{N},
    \label{eq-22}
\end{equation}
otherwise, it is calculated as 
\begin{equation}
    \rho_{j}^{\mathcal{M}} = \frac{12 \sum_{i=1}^{I} ( \sum_{i=1}^{I} R_{ijk} - \frac{1}{K} \sum_{i=1}^{I}\sum_{k=1}^{K} R_{ijk})^{2}}{I^{2}(K^{3}-K) - I\sum_{i=1}^{I} (T_{ij}^{3} - T_{ij}) }, \quad \forall j \in \mathcal{N},
    \label{eq-23}
\end{equation}
where $T_{ij}$ denotes the number of duplicate rankings from expert $i$ under attribute $j$.

According to Corollary \ref{corollary-05}, the weights assigned to attributes by experts are independent of alternative utilities. Thus, the correlation analysis for attributes only needs to consider experts, excluding alternatives. Kendall's W correlation coefficient for attributes among experts $\rho^{\mathcal{N}}$ can then be calculated in a similar manner. This coefficient ranges from 0 to 1, with values closer to 1 indicating higher correlation and greater consistency in expert decisions. Consequently, we propose the following test hypotheses for Kendall's W correlation coefficient:

\begin{itemize}
    \item[$-$] $H_0$: Inconsistent decision outcomes among experts;
    \item[$-$] $H_1$: Consistent decision outcomes among experts.
\end{itemize}

\begin{proposition}
    \citep{K76}.
    For small sample cases, $x_{j}^{\mathcal{M}} = \rho_{j}^{\mathcal{M}}(I-1)/(1-\rho_{j}^{\mathcal{M}}), \forall j \in 
    \mathcal{N}$ is approximately distributed as $F_{v_1,v_2}$, with the degrees of freedom $v_{1} = K - 1 - 2/I$ and $v_2 = (I-1)v_1$. The probability density function is given by 
    \[
        f(x;v_{1},v_{2})=\frac{\Gamma\left(\frac{v_1+v_2}{2}\right)\left(\frac{v_{1}}{v_{2}}\right)^{{\frac{v_{1}}{2}}}x^{{\frac{v_{1}}{2}-1}}}{\Gamma\left(\frac{v_1}{2}\right)\Gamma\left(\frac{v_2}{2}\right)\left(1+\frac{v_{1}}{v_{2}}x\right)^{{\frac{v_{1}+v_{2}}{2}}}}, \quad x>0.
    \]
    where $\Gamma$ is gamma function. $\rho^{\mathcal{N}}$ follows a similar approximation by replacing $K$ with $J$.
    \label{proposition-04}
\end{proposition}

By Proposition \ref{proposition-04}, for every $j \in \mathcal{N}$, if 
\[
    P\{ X > x_{j}^{\mathcal{M}} \} = 1 - F(x_{j}^{\mathcal{M}}; v_1, v_2) = 1 - \int_{0}^{x_{j}^{\mathcal{M}}} f(x; v_1, v_2) \mathrm{d} x \leq \alpha,
\]
then $H_0$ is rejected for alternatives under attribute $j$, indicating that decision outcomes for alternatives under attribute $j$ among experts are consistent. Kendall's W correlation coefficient for attributes among experts is tested similarly. Notably, failure to achieve consensus in expert hypothesis testing does not imply unacceptable decision outcomes. Divergent opinions among experts are common in group decision-making. Unlike pairwise comparison methods like AHP and BWM, GOPA avoids issues of consensus testing by deriving optimal weights from individual dominance ranking preference, mitigating the impact of conflicts among experts. The following gives the confidence level corresponding to Kendall’s W correlation.
\begin{definition}
    The global confidence level of GOPA is defined as
\begin{equation}
    GCL = LCL^{\mathcal{N}} \left( \sum_{j=1}^{J} w_{j}^{\mathcal{N}}LCL_j^{\mathcal{M}} \right),
    \label{eq-24}
\end{equation}
where $LCL^{\mathcal{N}}$ and $LCL_{j}^{\mathcal{M}}$ denote the local confidence levels (LCL) of attributes and alternatives under attributes, given by $LCL^{\mathcal{N}} = F(x^{\mathcal{N}}; v_1, v_2) \text{ and } LCL_j^{\mathcal{M}} = F(x_j^{\mathcal{M}}; v_1, v_2), \forall j \in \mathcal{N}$.
\end{definition}

The reason for adopting multiplicative aggregation of GCL is its ability to reveal a crucial characteristic in decision-making: any significant divergence in the decision outcomes will decrease overall decision reliability. The proposed confidence levels consist of two types: GCL and LCL. GCL assesses the overall reliability of GOPA decision outcomes, whereas the LCL evaluates the reliability of individual attributes and alternative outcomes. High sensitive problem occurs when the confidence level exceeds 0.99, very sensitive between [0.95, 0.99), sensitive between [0.90, 0.95), and less sensitive when below 0.9 \citep{M.J22}.

\section{Numerical Study}
\label{sec-4}
\subsection{Case Description and Data Collection}
\label{sec-4-1}

This study selects the improvisational emergency supplier selection (IESS) during the 7.20 mega-rainstorm disaster in Zhengzhou, China, as a case study to demonstrate the proposed GOPA \citep{P.Z22}. Unlike regular emergency supplier selection, IESS has distinct characteristics. Regular emergency supplier selection typically occurs during the emergency preparedness phase; however, the complexity and unpredictability of disasters may render pre-selected suppliers inadequate for disaster response. In such situations, IESS becomes essential. IESS requires rapid decision-making under tight time constraints and high uncertainty, making it challenging to obtain high-quality data and often necessitating reliance on partial preference information. Additionally, due to time pressure and information uncertainty, decision-makers' risk preferences significantly influence decision-making. Overall, the reason for selecting IESS to illustrate the applicability and effectiveness of GOPA is that it aligns with the GOPA decision analysis framework, and the urgency and complexity of decision-making in this scenario provide an ideal testing platform.

There are 10 emergency suppliers (labeled A1 to A10) available for selection within the disaster area, each with distinct characteristics. Some emphasize stable supply chains and prompt response capabilities despite higher costs, while others utilize strategic geographic positioning and efficient transport connections to expedite flood relief efforts during crises. Attributes selected for IESS include response speed (C1), delivery reliability (C2), geographic coverage (C3), operational sustainability (C4), collaborative experience and credibility (C5), and supply cost (C6) \citep{W.S.C.S.C.G24}. Five decision-makers (labeled E1 to E5) from various departments act as stakeholders for IESS, prioritized based on their decision-making authority as E5$>$E2$>$E1$>$E3$>$E4. These decision-makers provides their risk preference to each attribute, rankings of attributes and alternatives, and sufficiently determined partial preference information, as detailed in Online Supplemental Material. Table \ref{tab-C1} displays the frequency of occurrences of the expert global utility structure across each attribute. As discussed in Section \ref{sec-3-2-3}, the global utility structure in discrete prospects is characterized as risk-seeking (see Remark \ref{remark-01}). The uniform distribution is considered risk-neutral, while CARA and HARA are classified as risk-averse, and the S-shape represents a blend of risk-seeking and risk-averse. As for the parameters of the risk preference utility function in the continuous prospects, we use parameters consistent with Example \ref{example-03}.
\begin{table}[h]
\Scentering
\caption{Number of global utility structure occurrence under attributes}
\label{tab-C1}
\begin{tabular}{@{}llllllllll@{}}
\toprule
\multirow{2}{*}{Global   utility structure} & \multicolumn{5}{l}{Discrete prospect} & \multicolumn{4}{l}{Continuous prospect} \\ \cmidrule(l){2-10} 
      & RS & REF & RR & SR & ROC & Uniform & CARA & HARA & S-shape \\ \midrule
C1    & 0  & 0   & 0  & 0  & 1   & 1       & 0    & 1    & 2       \\
C2    & 0  & 1   & 0  & 0  & 2   & 1       & 0    & 1    & 0       \\
C3    & 0  & 0   & 1  & 0  & 1   & 0       & 3    & 0    & 0       \\
C4    & 1  & 0   & 1  & 1  & 0   & 0       & 1    & 1    & 0       \\
C5    & 2  & 0   & 0  & 0  & 1   & 0       & 1    & 1    & 0       \\
C6    & 1  & 0   & 0  & 2  & 2   & 0       & 0    & 0    & 0       \\
Total & 4  & 1   & 2  & 3  & 7   & 2       & 5    & 4    & 2       \\ \bottomrule
\end{tabular}
\end{table}

\subsection{Experiment Results}
\label{sec-4-2}

\begin{figure}[h]
\centering
\includegraphics[width=\textwidth]{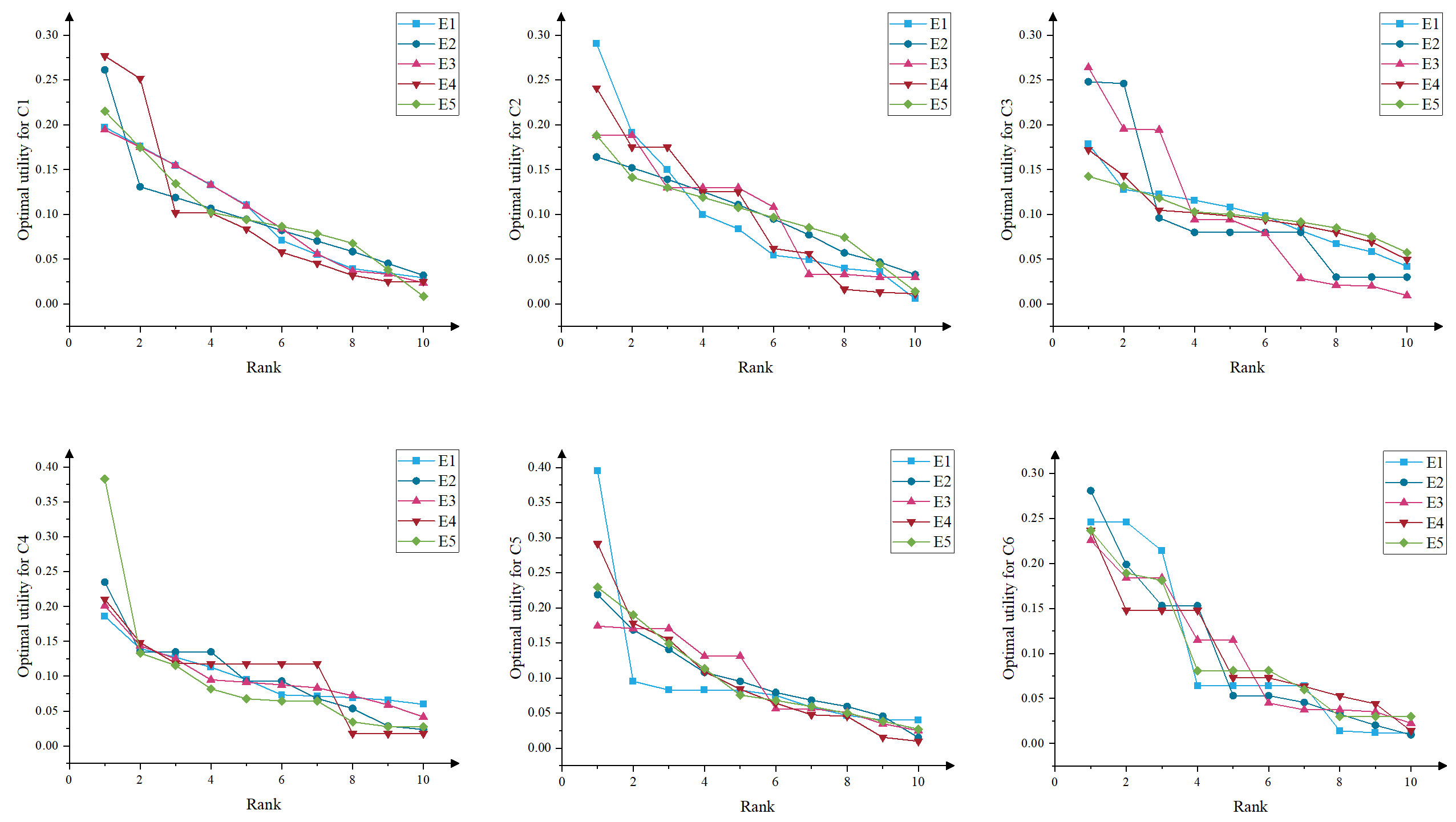}
\caption{Optimal utility of alternatives under different rankings under experts and attributes}
\label{fig-06}
\end{figure}

Figure \ref{fig-06} presents the optimal utilities of alternatives ranked across diverse attributes using partial preference information from experts. These results indicate that utility disparities among ranked alternatives are variable, underscoring the importance of partial preference information. Moreover, utilities derived from experts' partial preference information demonstrate non-convex and non-concave features, challenging the assumption that decision-makers' utilities adhere to specific functional structures in other utility elicitation research. Table \ref{tab-C2} presents the consistency measures of the expert group outcomes in decision-making. The IESS case overall is categorized as a less sensitive problem, with a GCL of 0.5797, despite high consistency among experts on attributes, with a LCL of 0.9951. However, notable disparities arise among experts in alternative outcomes under operational sustainability (C4). Therefore, initiating discussions of alternative assessments under operational sustainability (C4) is advisable.

\begin{table}[h]
\centering
\caption{Group consistency measures}
\label{tab-C2}
\begin{tabular}{@{}lllll@{}}
\toprule
ID & PSD    & \begin{tabular}[c]{@{}l@{}}Kendall's W \\ correlation coefficient\end{tabular} & LCL                      & GCL                      \\ \midrule
C1 & 0.0959 & 0.1893                                                                         & 0.4941                   & \multirow{16}{*}{0.5797} \\
C2  & 0.2089 & 0.2213 & 0.6355 &  \\
C3  & 0.0784 & 0.1496 & 0.3045 &  \\
C4  & 0.1595 & 0.0982 & 0.0993 &  \\
C5  & 0.3302 & 0.2882 & 0.8489 &  \\
C6  & 0.0660 & 0.2960 & 0.8658 &  \\
A1 & 0.1176 & \multirow{10}{*}{0.5154}                                                       & \multirow{10}{*}{0.9951} &                          \\
A2  & 0.1723 &        &        &  \\
A3  & 0.1345 &        &        &  \\
A4  & 0.1279 &        &        &  \\
A5  & 0.1207 &        &        &  \\
A6  & 0.0744 &        &        &  \\
A7  & 0.1862 &        &        &  \\
A8  & 0.1954 &        &        &  \\
A9  & 0.1773 &        &        &  \\
A10 & 0.1261 &        &        &  \\ \bottomrule
\end{tabular}
\end{table}

\begin{figure}[h]
\centering
\includegraphics[width=\textwidth]{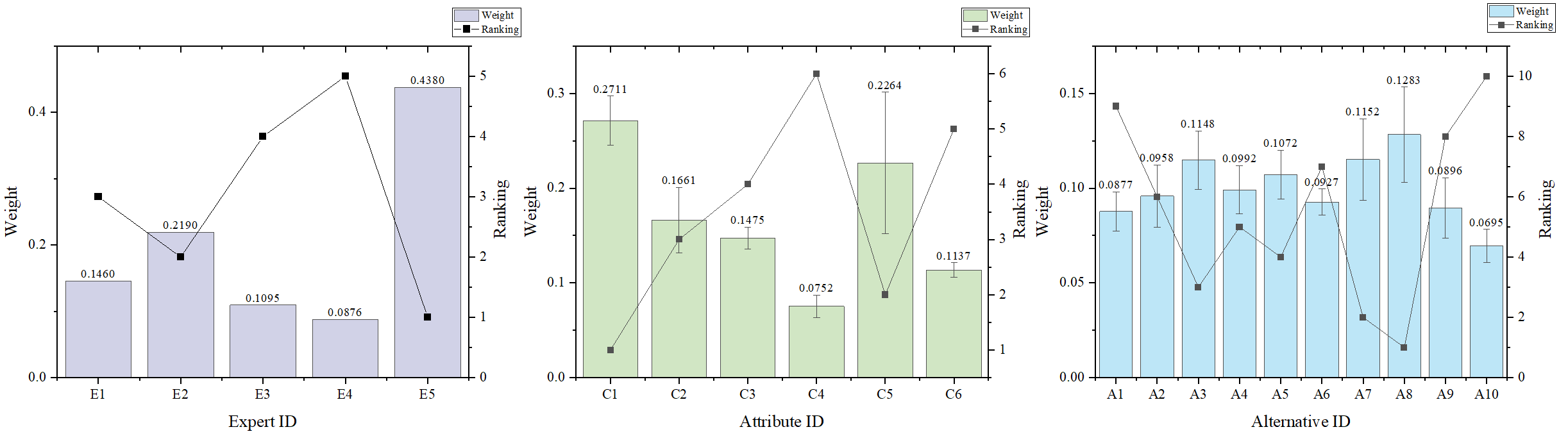}
\caption{Optimal weights of experts, attributes, and alternatives}
\label{fig-07}
\end{figure}

 Figure \ref{fig-07} depicts the weights assigned to experts, attributes, and alternatives. The expert weight results reveal a distinct hierarchy: E5 possesses the highest weight of 0.4380, succeeded by E2 at 0.2190. Following them, E1, E3, and E4 carry weights of 0.1460, 0.1095, and 0.0876, respectively. This distribution highlights a discernible pattern of diminishing marginal effects of RR weights across the experts. Regarding attributes, response speed (C1) emerges as the most critical factor, receiving a weight of 0.2444. This is closely followed by collaborative experience and credibility (C5) at 0.2007, delivery reliability (C2) at 0.1481, and geographic coverage (C3) at 0.1329. In contrast, supply cost (C6) and operational sustainability (C4) carry lower weights of 0.0931 and 0.0512, respectively. PSD of attributes indicates that only C5, at 33.02$\%$, significantly exceeds the threshold of 20$\%$, while other attributes exhibit normal levels of dispersion. It is noteworthy that, despite the increasing emphasis on aligning humanitarian operations with the UN Sustainable Development Goals, stakeholders in the case study perceive the sustainability of emergency suppliers as relatively insignificant. Regarding alternative weights, the top three alternatives rank as follows: A8, A7, and A3. Specifically, A8 holds the highest weight at 0.1283, followed by A7 at 0.1152 and A3 at 0.1148. PSDs for the alternatives are all below the 20$\%$ threshold, indicating stability. A8 excels in response speed (C1), operational sustainability (C4), and collaborative experience and credibility (C5), forming the basis for its selection. A7 complements A8 by demonstrating strengths in delivery reliability (C2), geographic coverage (C3), and supply cost (C6). Conversely, A3, overshadowed by A7 and A8, is not the preferred choice. In conclusion, A8 is identified as the optimal emergency supplier, with A7 as the secondary option.

\subsection{Sensitivity Analysis}
\label{sec-4-3}

\begin{figure}[h]
\centering
\includegraphics[width=\textwidth]{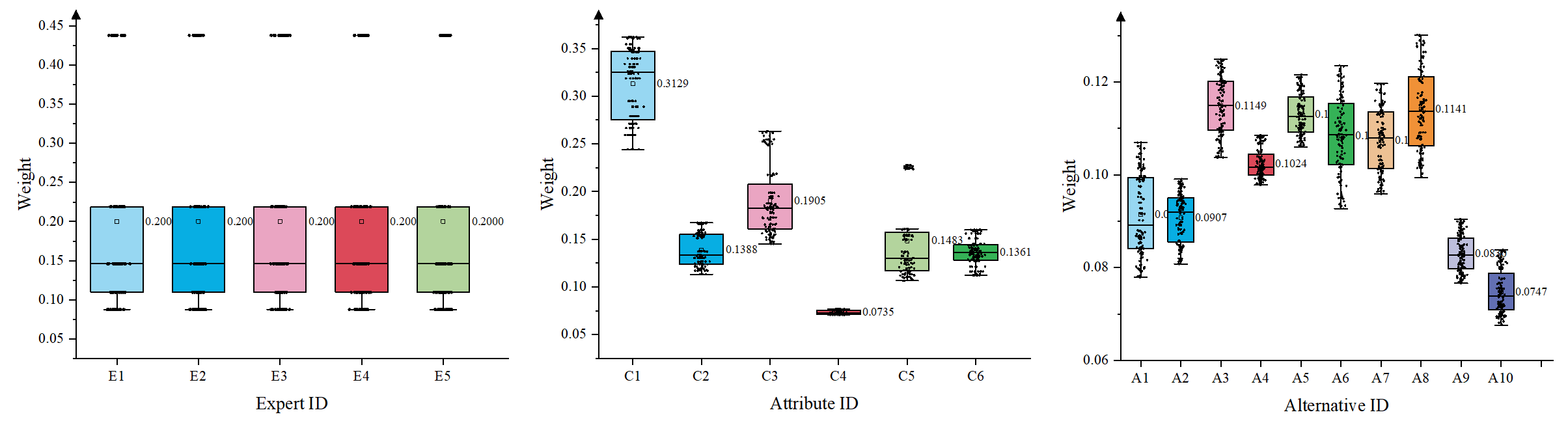}
\caption{Box plots of the weight outcomes}
\label{fig-08}
\end{figure}

Sensitivity analysis of input data or parameters is a vital numerical technique for evaluating the efficacy of MADM methods. This study specifically examines the sensitivity of expert rankings, particularly in the context of prioritizing five experts engaged in IESS. A permutation approach is utilized to create 120 potential scenarios to investigate the impact of varying rankings on decision outcomes. No interventions are made in the rankings of attributes and alternatives under each attribute given by experts, as these reflect their genuine preferences. Figure \ref{fig-08} and Table \ref{tab-C3} present detailed weight results and their descriptive statistics indicators.

\begin{table}[h]
\centering
\caption{Descriptive statistics of the weight outcomes}
\label{tab-C3}
\begin{tabular}{@{}lllllll@{}}
\toprule
    & Mean   & Skewness & Kurtosis & Coefficient of variation & Min    & Max    \\ \midrule
E1  & 0.2000 & 1.1019   & -0.3233  & 0.6380                   & 0.0876 & 0.4380 \\
E2  & 0.2000 & 1.1019   & -0.3233  & 0.6380                   & 0.0876 & 0.4380 \\
E3  & 0.2000 & 1.1019   & -0.3233  & 0.6380                   & 0.0876 & 0.4380 \\
E4  & 0.2000 & 1.1019   & -0.3233  & 0.6380                   & 0.0876 & 0.4380 \\
E5  & 0.2000 & 1.1019   & -0.3233  & 0.6380                   & 0.0876 & 0.4380 \\
C1  & 0.3129 & -0.3725  & -1.3224  & 0.1209                   & 0.2443 & 0.3620 \\
C2  & 0.1388 & 0.2554   & -1.4153  & 0.1238                   & 0.1132 & 0.1673 \\
C3  & 0.1905 & 0.8466   & -0.6099  & 0.1927                   & 0.1450 & 0.2634 \\
C4  & 0.0735 & 0.2999   & -1.4365  & 0.0289                   & 0.0707 & 0.0770 \\
C5  & 0.1483 & 1.0830   & -0.3420  & 0.2805                   & 0.1065 & 0.2279 \\
C6  & 0.1361 & 0.0000   & -0.7241  & 0.1009                   & 0.1122 & 0.1599 \\
A1  & 0.0914 & 0.1745   & -1.3541  & 0.0956                   & 0.0779 & 0.1070 \\
A2  & 0.0907 & -0.2485  & -1.2175  & 0.0586                   & 0.0808 & 0.0991 \\
A3  & 0.1149 & -0.1332  & -1.1672  & 0.0536                   & 0.1037 & 0.1249 \\
A4  & 0.1024 & 0.5820   & -0.8347  & 0.0305                   & 0.0978 & 0.1085 \\
A5  & 0.1130 & 0.2920   & -1.0886  & 0.0400                   & 0.1060 & 0.1216 \\
A6  & 0.1084 & -0.0949  & -0.9048  & 0.0786                   & 0.0927 & 0.1235 \\
A7  & 0.1074 & 0.0122   & -1.2450  & 0.0649                   & 0.0959 & 0.1197 \\
A8  & 0.1141 & 0.1595   & -1.1003  & 0.0768                   & 0.0995 & 0.1302 \\
A9  & 0.0830 & 0.2469   & -1.1072  & 0.0476                   & 0.0766 & 0.0904 \\
A10 & 0.0747 & 0.4826   & -1.0180  & 0.0628                   & 0.0676 & 0.0838 \\ \bottomrule
\end{tabular}
\end{table}

Descriptive statistical analysis of expert weights reveals significant consistency among the averages of five experts, averaging 0.2. The range of expert weights spans from a minimum of 0.0876 to a maximum of 0.4380, indicating their consistency and stability in weight allocation. This consistency is further validated across all other assessment metrics, aligning with observations from the entire ranking experimental design. Notably, the equal frequency of appearances by each expert across different ranking positions further reinforces the consistency of descriptive statistical results among the experts. Among the attributes, C1 has the highest average at 0.3129, followed by C3, C5, C2, and C6 at 0.1905, 0.1483, 0.1388, and 0.1361, respectively. C4 exhibits the lowest average at 0.0735. Regarding skewness, most attributes demonstrate a right-skewed tendency, indicating data primarily distributed on the positive side, although C6 approaches symmetry, while C1 shows a left-skewed distribution. C3, C5, and C6 exhibit higher kurtosis, while C1, C2, and C4 are relatively flat. Among these attributes, only C5's coefficient of variation is 0.2805, significantly exceeding the threshold of 0.2, while the others demonstrate comparatively stable dispersion. Based on the average weights of the alternatives, the top four rankings are A3 (0.1149), A8 (0.1141), A5 (0.1130), and A6 (0.1084). The coefficient of variation for all alternatives does not exceed the specified threshold. Moreover, the weights of A1, A4, A5, A7, A8, and A9 show positive skewness, indicating a bias towards higher values. Conversely, the remaining alternatives display negative skewness. The kurtosis of the alternatives is significantly lower compared to that of the attributes, suggesting a relatively flat distribution of weights among the alternatives. In conclusion, the weight outcomes for GOPA of sensitivity analysis are relatively stable.

\subsection{Comparison Analysis}
\label{sec-4-4}

This study validates GOPA by comparing its ranking outcomes with those of 7 other methods: OPA, CODAS, COPRAS, MACBETH, MAIRCA, MARCOS, and TOPSIS. These methods are chosen for their prominence as baseline approaches in the current MADM field. Except for OPA, all methods use a decision matrix of evaluation scores and require predefined attribute weights. According to Proposition \ref{proposition-02}, the weights allocated to alternatives across different attributes $w_{jk}$ can be decomposed into attribute weights $w_{j}^{\mathcal{N}}$ and alternative utilities based on partial preference information $u_{jk}$, where $u_{jk}$ can be regarded as the unweighted utility. Thus, $w_{j}^{\mathcal{N}}$ and $u_{jk}$ serve as the predefined weights and decision data for these methods, considering the risk preference and partial preference information. Data recollection is avoided to prevent increasing the cognitive burden and uncertainty for experts who have already provided sufficiently determined preference information. The input data for the comparison methods can be found in Supplemental Material. Table \ref{tab-C4} presents the ranking results of multiple methods.

\begin{table}[h]
\Scentering
\renewcommand{\arraystretch}{0.75}
\caption{Alternative rankings of multiple MADM methods}
\label{tab-C4}
\begin{tabular}{@{}lllllllll@{}}
\toprule
    & GOPA & OPA & CODAS & COPRAS & MACBETH & MAIRCA & MARCOS & TOPSIS \\ \midrule
A1  & 9    & 9   & 9     & 9      & 9       & 9      & 9      & 9      \\
A2  & 6    & 5   & 5     & 6      & 7       & 7      & 7      & 5      \\
A3  & 3    & 3   & 2     & 3      & 1       & 1      & 2      & 2      \\
A4  & 5    & 6   & 7     & 5      & 5       & 5      & 5      & 6      \\
A5  & 4    & 2   & 4     & 4      & 4       & 4      & 4      & 4      \\
A6  & 7    & 8   & 6     & 7      & 6       & 6      & 6      & 8      \\
A7  & 2    & 4   & 3     & 2      & 2       & 2      & 3      & 3      \\
A8  & 1    & 1   & 1     & 1      & 3       & 3      & 1      & 1      \\
A9  & 8    & 7   & 8     & 8      & 8       & 8      & 8      & 7      \\
A10 & 10   & 10  & 10    & 10     & 10      & 10     & 10     & 10     \\ \bottomrule
\end{tabular}
\end{table}

\begin{figure}[h]
\centering
\includegraphics[width=3.75in]{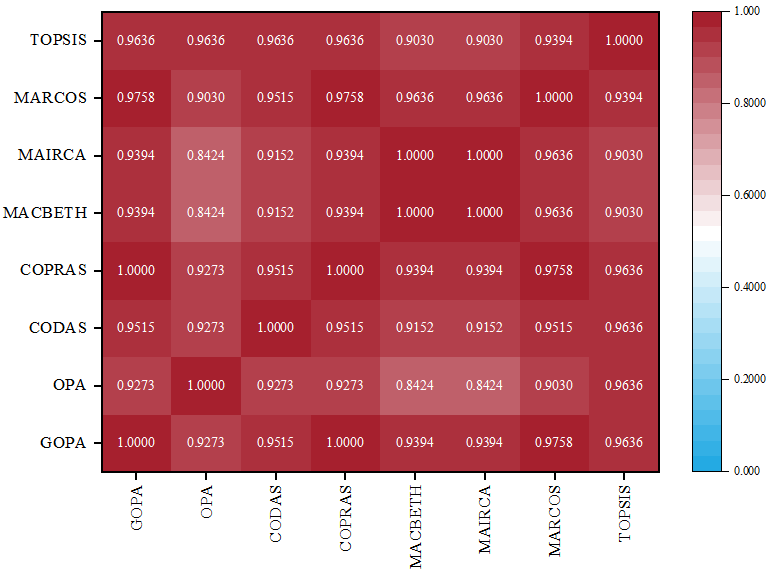}
\caption{Spearman correlation heat map of alternative rankings of multiple MADM methods}
\label{fig-09}
\end{figure}

This study then utilizes Spearman correlation coefficient to evaluate the correlation between these ordinal sequences of alternative rankings \citep{W24}. As the value approaches 1, the two rankings indicate a strong positive correlation. Conversely, when the value approaches -1, it signifies a significant negative correlation between the rankings. The value close to 0 suggests a lack of apparent correlation between the rankings. Figure \ref{fig-09} shows the Spearman correlation heat map of multiple MADM methods. The analysis reveals significant positive correlations among the eight MCDA methods, substantiating their similarity and potential complementarity. GOPA shows high correlation coefficients above 0.9 with other methods, significant at the 1$\%$ level. In contrast, OPA demonstrates lower correlations due to its neglect of decision makers' risk preferences and partial preference information. Other methods utilize GOPA-derived outputs to generate diverse rankings, typically requiring complete attribute information and predefined weights. However, reliance on alternative rankings per attribute, partial preference data, and expert risk preferences may compromise decision reliability. In summary, the GOPA method facilitates dependable decisions with minimal data requirements.

\section{Extension and Discussion}
\label{sec-5}

\subsection{Stochastic Dominance Constraints in Weight Elicitation}
\label{sec-5-1}

This section extends subjective preference-based weight elicitation to optimization with stochastic dominance constraints. Over the past few decades, stochastic dominance constraint optimization has emerged as a critical method for managing risk, particularly in the extensive application of portfolio selection problems \citep{P.D24}. The classic approach maximizes a return function while ensuring that a controlled stochastic performance function stochastically dominates a specified reference benchmark \citep{D.R03}. As noted in Remark \ref{remark-03}, the analytical solution of GOPA provides a reference benchmark for subjective preference-based weight elicitation, forming the foundation for transforming it into an optimization problem with stochastic dominance constraints. We consider the expert opinions as the realization with corresponding probability $p_{i} = W^{\mathcal{Q}}_{i}$. Then, we introduce the following optimization with stochastic dominance constraints:
\begin{equation}
        \max\limits_{\mathbf{w} \in \mathcal{W}} f(\mathbf{w}), \quad \text{s.t. } z^{*} \preceq_{(2)} g(\mathbf{w}),
        \label{eq-25}
\end{equation}
where $g(w_{jr}) = \left(s_{j}w_{jr}\right) \big/ \left(K_{j} U^{*}_{jr}\right)$, $f: \mathcal{W} \mapsto \mathbb{R}$ is a concave objective function, $z^{*}$ and $U^{*}_{jr}$ are constants derived from analytical solution in Theorem \ref{theorem-05}. 

\begin{lemma}
    The second order stochastic dominance $z^{*} \preceq_{(2)} g(\mathbf{w})$ is equivalent to

    1) $\mathbb{E}[(\eta - z^{*})^{+})] \ge \mathbb{E}[(\eta - g(\mathbf{w}))^{+})], \forall \eta \in \mathbb{R}$;

    2) $\mathbb{E}[u(z^{*})] \leq \mathbb{E}[u(g(\mathbf{w}))]$ for all concave non-decreasing utility functions $u: \mathbb{R} \mapsto \mathbb{R}$ such that both expectations exists.
    \label{lemma-02}
\end{lemma}

The second equivalence ensures that the weight assignments dominate the reference benchmark for all risk-averse utility functions, which satisfies the weight discrimination requirements. When consider the multi-attribute individual decision-making problem, the second order dominance reduces to the first order dominance. In addition, the objective function $f$ can be given by the common objective function employed in approximate weighting methods for MADM, such as entropy, squared deviation, and absolute dominance degree \citep{A17, A24}. By the first equivalence in Lemma \ref{lemma-02}, we have 
\begin{equation}
        \max\limits_{\mathbf{w} \in \mathcal{W}} f(\mathbf{w}), \quad \text{s.t. } \mathbb{E}[u(z^{*})] \leq \mathbb{E}[u(g(\mathbf{w}))].
        \label{eq-26}
\end{equation}
Note that $z^{*}$ is not a random variable with only one realization, which yields $\mathbb{E}[u(z^{*})] =u(z^{*})$. Consider the scenario where the actual expert weight distribution $P^{*}\in \mathcal{P}$ is uncertain and may differ from the empirical distribution $\Hat{P}$ (i.e., $W^{\mathcal{Q}}$). By the spirit of emerging robust satisficing model \citep{L.S.Z23}, we can sacrifice some performance (weight differentiation) to enhance robustness against uncertainty among expert weights, expressed as:
\begin{equation}
    \max\limits_{\mathbf{w} \in \mathcal{W}, \theta} \theta, \quad \text{s.t. } \alpha u(z^{*}) - \mathbb{E}[u(g(\mathbf{w}))] \ge \theta \Delta(\mathbb{P}, \Hat{\mathbb{P}}), \quad \forall \mathbb{P} \in \mathcal{P}.
\end{equation}
where $\alpha \in (0,1] $ represents the tolerance level. The above discussion on optimization with stochastic dominance constraints and robust satisficing framework offers a potential approach to address large-scale multi-attribute group decision-making problem, which we leave for future exploration.

\subsection{Accounting for Weight Elicitation Errors}
\label{sec-5-2}

This section discuss the potential elicitation errors in decision weights within GOPA. These errors may stem from various causes, such as the decision-maker’s preferences not aligning with the axioms of expected utility theory or contaminated preference information \citep{A.D15}. Given the critical nature of these axioms of expected utility theory and the existence of imprecise information, GOPA is prescriptive. Specifically, our main objective is to aid decision-makers who believe in expected utility theory and have confidence in their preference information to determine decision weights that accurately reflect their preferences. Similar to the spirit of \citet{B.O13}, when a decision-maker’s preferences do not satisfy a particular axiom, GOPA should elicit and adopt decision weights that best explain these inconsistencies. Practically, these inconsistencies can be viewed as minor measurement errors needing correction to accurately identify the decision-maker’s true intentions, even if they do not explicitly articulate their preferences. When inconsistencies are detected, we suggest considering error margins in the formula. We provide three interpretable types of errors and corresponding corrective formulations.

When considering the noise of the differentiation of the decision weights, $g(\mathbf{w})$ is evaluated as $g(\mathbf{w}) + \gamma$ for a perturbation $\gamma \in \mathbb{R}$. Then, we can relax the constraint $f(z, \mathbf{u}^{*}_{ijr}) \leq g(\mathbf{w}_{ijr})$ with
\begin{equation}
    f(z, \mathbf{u}^{*}_{ijr}) - g(\mathbf{w}_{ijr}) \leq \gamma_{ijr}, \quad \forall (i,j,r) \in \mathcal{U},
    \label{eq-27}
\end{equation}

where $\gamma_{ijr} \ge 0$ is the error tolerance for the alternative ranked $r$ on attribute $j$ under the preference of expert $i$. Given a budget $\Gamma \ge 0$, the total error $\sum_{i=1}^{I}\sum_{j=1}^{J}\sum_{r=1}^{K_{ij}} \gamma_{ijr} \leq \Gamma$ is considered as a measure of inconsistency for the feasibility of the second-stage GOPA problem. The objective of Equation (\ref{eq-14}) is then replaced by 
\begin{equation}
    \underset{z, \mathbf{w}, \mathbf{\gamma}}{\mathrm{max}} \text{ } z + \sum_{i=1}^{I}\sum_{j=1}^{J}\sum_{r=1}^{K_{ij}} \gamma_{ijr}.
    \label{eq-28}
\end{equation}

When accounting for errors in utility elicitation for alternatives in the first-stage problem of GOPA, for continuous prospects, the variable $r$ is adjusted to $r + \tau_{ijr}$, where $\tau_{ijr} \in (-1,1)$ is a given parameter. For discrete prospects, interpolation is used to elicit utility. Consequently, the constraint $f(z, \mathbf{u}^{*}_{ijr}) \leq g(\mathbf{w}_{ijr})$ for the specific $\Bar{r}$ is replaced with $f(z, \mathbf{u}^{\Bar{*}}_{ijr}) \leq g(\mathbf{w}_{ijr})$, where
\[
    U_{ijr}^{\Bar{*}}= \begin{cases} U_{ij}^{*}(r + \tau_{ijr}), & \forall(i,j,r) \in \mathcal{U}^1, \\
    \left( \int_{0}^{r + \tau_{ijr}}u_{ij}^{*}(K_{ij}-x)dx \right) \Big/ \left( \sum_{r=1}^{K_{ij}}\int_{0}^{r + \tau_{ijr}}u_{ij}^{*}(K_{ij}-x)dx \right),  & \forall(i,j,r) \in \mathcal{U}^2.
    \end{cases}
\]

A third option posits that $1 - \epsilon_{ij}$ of the $K_{ij}$ dominance relations are correct, meaning expert $i$ errs in at most $\epsilon_{ij} K_{ij}$ relations on attribute $j$. We can introduce a binary variable $\vartheta_{ijr}$, which equals 1 if expert $i$ is incorrect about the relation of ranking $r$ on attribute $j$, subject to the constraint $\sum_{r = 1}^{K_{ij}} \vartheta_{ijr} \leq \epsilon_{ij} K_{ij}, \forall (i,j) \in \mathcal{U}$. The constraint $f(z, \mathbf{u}^{*}_{ijr}) \leq g(\mathbf{w}_{ijr})$ is then replaced by two constraints: $(1 - \vartheta_{ijr}) M + f(z, \mathbf{u}^{*}_{ijr}) \ge g(\mathbf{w}_{ijr})$ and $f(z, \mathbf{u}^{*}_{ijr}) \leq \vartheta_{ijr} M + g(\mathbf{w}_{ijr}), \forall (i,j,r) \in \mathcal{U}$, where $M$ is a large constant. This second-stage problem of GOPA can be solved using the cutting-plane method.

\subsection{Advantage and Insight of Generalized Ordinal Priority Approach}
\label{sec-5-3}

Currently, MADM methods are broadly categorized into two groups: weighting methods and ranking methods. Weighting methods like AHP, FUCOM, LBWA, and BWM typically assign separate weights to experts, attributes, and alternatives. Ranking methods, such as TOPSIS, VIKOR, PROMETHEE, COPRAS, EDAS, and MABAC, generally prioritize alternatives based on their performance and pre-obtained weights of experts and attributes. Consequently, many existing MADM studies integrate both weighting and ranking methods. In contrast, GOPA represents a comprehensive MADM approach that determines weights for experts, attributes, and alternatives simultaneously for prioritizing alternatives, without the need for expert opinion aggregation or pre-acquired weights.

Most MADM methods rely on complete preference information, as discussed in Section \ref{sec-3-2-1} for decision-making. For example, pairwise comparison methods such as AHP, BWM, and FUCOM are typically based on ratio scale preferences, while quasi-distance-based methods like TOPSIS, VIKOR, EDAS, and COPRAS generally use absolute difference preferences \citep{Z.K.Y.Z.D23}. Furthermore, fitting utility functions based on expected utility theory is often achieved by pairwise comparison of a large number of lotteries under parameter assumptions \citep{G.X.Z23}. However, decision-making based on complete preference information is often idealized and expensive, which faces the challenges of time constraints, insufficient data, and expert cognitive burden. While some MADM studies incorporate fuzzy theory and grey theory to address uncertainty, they do not eliminate the data requirement \cite{W24}. Moreover, methods proposed by \citet{A15, A17, A24} combine extremal points of weight constraint derived from partial preference information to determine optimal decision weights but do not address MADM involving multiple experts and attributes with risk preferences. In some sense, GOPA also combines weights from partial preference information through the “estimate-then-optimize” contextual optimization framework. However, unlike other methods, GOPA comprehensively considers risk preferences across different experts and attributes. Thus, GOPA could inspire a kind of two-stage optimization-based MADM method with a general principle: eliciting the utility of alternatives with the risk preference and partial preference information of decision-makers across attributes in the first stage, then optimizing weights in a normalized weight space including experts, attributes, and alternatives. The IESS case aptly clarifies the practical conditions for applying GOPA, characterized by time and cost constraints in decision-making and pronounced risk preferences among decision-makers. Overall, GOPA shows the following advantages:
\begin{itemize}
    \item[1)] It utilizes easily obtainable and stable ranking data and sufficiently determined partial preference information as inputs for the model;
    \item[2)] It personalizes the decision weights based on partial preference information of various attributes to accurately reflect decision-makers' risk preferences in situations with incomplete information;
    \item[3)] It concurrently determines the weights of experts, attributes, and alternatives, eliminating the necessity for data normalization, expert opinion aggregation, and weight pre-acquisition techniques.
\end{itemize}

\section{Conclusion}
\label{sec-6}

In this study, we restate OPA, an emerging MADM method, and derive its properties, including solution efficiency, analytical solution expression and its decomposability, and the relationship between optimal decision weights and rank-based surrogate weights. It is notable that ROC weights are found in the reformulated modeling of OPA. Building upon this, we extend ROC weights within OPA to a more general utility structure, proposing GOPA. GOPA employs an “estimate-optimize” contextual optimization framework to derive decision weights for experts, attributes, and alternatives under incomplete information. In the first stage of GOPA, we explore utility elicitation based on general partial preference information and risk attitudes under discrete and continuous prospects. To manage partial preference information with limited data that do not meet statistical requirements, this study employs cross-entropy utility minimization to derive optimal utility, incorporating partial preference information of weak ordered relations, absolute differences, ratio scales, and lower bounds. Rank-based surrogate weights and risk preference utility functions, widely adopted in current decision science domains, serve as global target utility structures for discrete and continuous prospects, respectively. The utilities derived from various experts and attributes in the first stage serve as parameters for optimizing decision weights for experts, attributes, and alternatives in the second stage. We present validation metrics for GOPA group decisions, including the percentage standard deviation, correlation coefficients, and confidence level measurements. Leveraging the analytical solution of GOPA, we discuss the conversion of subjective preference-based utility elicitation problems into optimization with stochastic dominance constraints and robust satisficing, offering a promising avenue for addressing large-scale group decision-making challenges. Additionally, we discuss potential elicitation errors in GOPA and propose corresponding corrective formulations. The improvisational emergency supplier selection during the 7.20 mega-rainstorm disaster in Zhengzhou, China is presented as a case study to demonstrate and validate GOPA.

Theoretical analysis reveals several advantageous properties of GOPA:
\begin{itemize}
    \item [1)]Model generalizability: GOPA integrates mainstream utility elicitation approaches within the MADM domain. Through its modeling process, GOPA can degenerate into OPA, ranking-based alternative weights, and risk preference utility functions with appropriate parameter selection. When global utility structures are absent, minimizing cross-entropy utility transforms into maximizing entropy utility from partial preference information. In this sense, GOPA serves as an alternative modeling approach for utility elicitation based on partial preference information.
    \item [2)] Analytical solvability: The analytical solution of the second-stage problem in GOPA can be solved through the Lagrange multiplier method. It is observed that the optimal value of the objective function remains constant under various global utility structures and partial preference information, regardless of discrete or continuous prospects. This property forms the foundation for deriving expressions for other weight elicitation methods considering multi-attributes and multi-experts by substituting the utility distribution of alternatives into GOPA's analytical solution expression.
    \item [3)] Risk preference independence: The analytical solution of GOPA demonstrates that decision-makers' risk preferences do not influence the weight allocation for decision-makers and attributes, only the weight outcomes of alternatives. Additionally, GOPA's modeling process ensures that a decision-maker's risk preference for one attribute does not affect their preferences for other attributes or those of other decision-makers. Overall, in GOPA, decision-makers' risk preferences are independent and do not impact outcomes beyond alternative weight allocation.
\end{itemize}

We would like to highlight several limitations and corresponding future research directions in addition to the extensions of GOPA in prospect theory and local utility structures discussed in this study. First, this study does not explore how to elicit experts' global utility structures and partial preferences more effectively, concentrating on the model formulation and analysis. Insights from finance and behavioral economics could enhance application-specific integration of GOPA in the above limitation. Second, the current version of GOPA primarily considers subjective decision information, excluding objective decision data. Notably, the first-stage derived utility distribution for alternatives in GOPA can be considered a generalized standardized performance indicator. Thus, selecting appropriate data standardization techniques for objective decision data of diverse scales and dimensions, and theoretically exploring the relationship between GOPA's solution space and the decision values of quasi-distance-based MADM methods, presents a interesting research opportunity. This exploration could contribute to a more unified MADM framework. Finally, in large-scale group decision-making scenarios, extending GOPA to account for uncertain information and elicitation errors is also noteworthy. In this scenario, combining machine learning classification, contextual optimization, and distributionally robust optimization could offer a promising approach to GOPA modeling.

\section*{Author Contributions}
\textbf{Renlong Wang:} Writing-review $\&$ editing, Writing-original draft, Validation, Software, Methodology, Investigation, Conceptualization.

\section*{Declaration of competing interest}
The authors declare that they have no known competing financial interests or personal relationships that could have appeared to influence the work reported in this paper.

\section*{Data Availability}
Data will be available upon request or in the Supplemental Material on the GitHub software repository (\href{https://github.com/Renlong-WANG/GOPA}{https://github.com/Renlong-WANG/GOPA}).

\bibliographystyle{elsarticle-harv} 
\bibliography{GOPA}

\newpage

\appendix
\section{Proofs for Section \ref{sec-2}}

\newproof{lemma-proof-01}{Proof of Lemma \ref{lemma-01}}
\begin{lemma-proof-01}
By the definition of problem of the second-stage problem, we have $z^{*} \le \delta(\mathbf{w})$ for every $(i,j,k)\in {{\mathcal{X}}^{1}} \cup {{\mathcal{X}}^{2}}$. It follows that $z^{*} \le \min \delta(\mathbf{w})$. Assume that $\min \delta(\mathbf{w})$ is strictly bigger than $z^{*}$. However, $\mathbf{w} \in FS(P)$, and this will increase the objective value of the first-stage problem. Contradiction. 
\qed
\end{lemma-proof-01}

\newproof{theorem-proof-01}{Proof of Theorem \ref{theorem-01}}
\begin{theorem-proof-01}
Let $\mathbf{w} \in OS(Q)$. Assume that $\mathbf{w}$ is not efficient. Then, there exists $\mathbf{w}' \in FS(Q)$ such that $\delta(\mathbf{w}') > \delta(\mathbf{w})$ for some $(i,j,k)' \in \mathcal{Y}$ and $\delta(\mathbf{w}') \ge \delta(\mathbf{w})$ for every $(i,j,k) \in \mathcal{Y} \backslash \{(i,j,k)'\}$. Then, we have
\[
\sum\limits_{i \in \mathcal{Q}} \sum\limits_{j \in \mathcal{N}} \sum\limits_{k \in \mathcal{M}} \delta(\mathbf{w}') > \sum\limits_{i \in \mathcal{Q}} \sum\limits_{j \in \mathcal{N}} \sum\limits_{k \in \mathcal{M}} \delta(\mathbf{w}) = z^{*}
\]
which contradicts the assumption that $\mathbf{w} \in OS(Q)$. 
\qed
\end{theorem-proof-01}

\newproof{theorem-proof-02}{Proof of Theorem \ref{theorem-02}}
\begin{theorem-proof-02}
Consider the reformulated OPA model in Equation (\ref{eq-09}), which is a typical convex optimization problem with its feasible region being convex and linear objective function being both convex and concave. Therefore, we can employ the Lagrange multiplier method and obtain the following Lagrange function.
\begin{equation}
\begin{aligned}
L(z,{{w}_{ijr}},\alpha ,{{\beta }_{ijr}})=z - & \alpha \left(\sum\limits_{i=1}^{I}{\sum\limits_{j=1}^{J}{\sum\limits_{r=1}^{{{K}_{ij}}}{{{c}_{ijr}}{{w}_{ijr}}}}}-1 \right) \\
& -\sum\limits_{i=1}^{I}{\sum\limits_{j=1}^{J}{\sum\limits_{r=1}^{{{K}_{ij}}}{{{\beta }_{ijr}} \left({{t}_{i}}{{s}_{ij}}{{w}_{ijr}}-\left(\sum\nolimits_{h=r}^{{{K}_{ij}}}{\frac{1}{h}} \right)z\right)}}}  
\end{aligned}
\label{eq-sub32-04}
\end{equation}
Let 
\begin{equation}
\frac{\partial L(z,{{w}_{ijr}},\alpha ,{\beta }_{ijr})}{\partial \alpha }=\frac{\partial L(z,{{w}_{ijr}},\alpha ,{\beta }_{ijr})}{{{\beta }_{ijr}}}=0.
\label{eq-sub32-05}
\end{equation}

These yield
\begin{equation}
\frac{\partial L(z,{{w}_{ijr}},\alpha ,{\beta }_{ijr})}{\partial \alpha }=\sum\limits_{i=1}^{I}{\sum\limits_{j=1}^{J}{\sum\limits_{r=1}^{{{K}_{ij}}}{{{c}_{ijr}}{{w}_{ijr}}}}}-1=0
\label{eq-sub32-06}
\end{equation}
and
\begin{equation}
\frac{\partial L(z,{{w}_{ijr}},\alpha ,{\beta }_{ijr})}{{\beta }_{ijr}}={{t}_{i}}{{s}_{ij}}{{w}_{ijr}}-\left(\sum\nolimits_{h=r}^{{{K}_{ij}}}{\frac{1}{h}}\right)z=0, \quad \forall (i,j,r)\in \mathcal{U}.
\label{eq-sub32-07}
\end{equation}

Substituting Equation (\ref{eq-sub32-07}) into Equation (\ref{eq-sub32-06}) yields
\begin{equation}
\left(\sum\limits_{i=1}^{I}{\sum\limits_{j=1}^{J}{\sum\limits_{r=1}^{{{K}_{ij}}}{\frac{{{c}_{ijr}}\sum\nolimits_{h=r}^{{{K}_{ij}}}{\frac{1}{h}}}{{{t}_{i}}{{s}_{ij}}}}}}\right) z = 1 \Leftrightarrow z= 1\Bigg/ \sum\limits_{i=1}^{I}{\sum\limits_{j=1}^{J}{\sum\limits_{r=1}^{{{K}_{ij}}}{\frac{{{c}_{ijr}}\sum\nolimits_{h=r}^{{{K}_{ij}}}{\frac{1}{h}}}{{{t}_{i}}{{s}_{ij}}}}}} .
\label{eq-sub32-08}
\end{equation}

Then, substituting Equation (\ref{eq-sub32-08}) into Equation (\ref{eq-sub32-07}) yields
\begin{equation}
{{w}_{ijr}}=\frac{\sum\nolimits_{h=r}^{{{K}_{ij}}}{\frac{1}{h}}}{{{t}_{i}}{{s}_{ij}}}z = \left(\sum\nolimits_{h=r}^{{{K}_{ij}}}{\frac{1}{h}} \right) \Bigg/ \left({{t}_{i}}{{s}_{ij}}\sum\limits_{i=1}^{I}{\sum\limits_{j=1}^{J}{\sum\limits_{r=1}^{{{K}_{ij}}}{\frac{{{c}_{ijr}}\sum\nolimits_{h=r}^{{{K}_{ij}}}{\frac{1}{h}}}{{{t}_{i}}{{s}_{ij}}}}}} \right), \quad \forall (i,j,r)\in \mathcal{U}.
\label{eq-sub32-09}
\end{equation}

Finally, we have the analytical solution $z^{*}$ and $w^{*}_{ijk}, \forall (i,j,r) \in \mathcal{U}$ of OPA. \qed
\label{theorem-proof-02}
\end{theorem-proof-02}

\newproof{corollary-proof-01}{Proof of Corollary \ref{corollary-01}}
\begin{corollary-proof-01}
By the analytical solution of OPA with no missing and duplicate rankings, we have
\[
    {{w}^{*}_{ijr}}=\frac{\frac{1}{K}\sum\nolimits_{h=r}^{K}{\frac{1}{h}}}{{{t}_{i}}{{s}_{ij}} \left(\sum\nolimits_{p=1}^{I}{\frac{1}{p}}\right)\left(\sum\nolimits_{q=1}^{J}{\frac{1}{q}}\right)}=\frac{v_{r}^{ROC}}{{{t}_{i}}{{s}_{ij}}\left(\sum\nolimits_{p=1}^{I}{\frac{1}{p}}\right)\left(\sum\nolimits_{q=1}^{J}{\frac{1}{q}}\right)}, \quad \forall (i,j,r)\in \mathcal{U}.
\]
When there is only one expert and attribute or when the importance of all experts and attributes is equal, the alternative weights in OPA will reduce to the rank order centroid weights.
\qed
\end{corollary-proof-01}

\newproof{corollary-proof-02}{Proof of Corollary \ref{corollary-02}}
\begin{corollary-proof-02}
Consider the case without missing and duplicate rankings. By summing the alternatives, we have the optimal weights of attributes under experts:
\[
    {{w}_{is}^{*}}=\sum\limits_{r=1}^{K}{\frac{\frac{1}{K}\sum\nolimits_{h=r}^{K}{\frac{1}{h}}}{{{t}_{i}}{{s}_{ij}}\left(\sum\nolimits_{p=1}^{I}{\frac{1}{p}}\right)\left(\sum\nolimits_{q=1}^{J}{\frac{1}{q}}\right)}}=\frac{1}{{{t}_{i}}{{s}_{ij}}\left(\sum\nolimits_{p=1}^{I}{\frac{1}{p}}\right)\left(\sum\nolimits_{q=1}^{J}{\frac{1}{q}}\right)}=\frac{v_{s}^{RR}}{{{t}_{i}} \left(\sum\nolimits_{p=1}^{I}{\frac{1}{p}}\right)}, \quad (i,s) \in \mathcal{U}.
\]
When there is only one expert or when the importance of all experts is equal, the attribute weights derived from OPA will reduce to the rank reciprocal weights. 
\qed
\end{corollary-proof-02}

\newproof{proposition-proof-02}{Proof of Proposition \ref{proposition-02}}
\begin{proposition-proof-02}
Consider the case without missing and duplicate rankings. By summing the attributes, we have the weights of alternatives under experts:
\[
    w_{ir}^{*} = \sum_{j=1}^J\frac{\sum_{h=r}^K\frac{1}{h}}{t_is_{ij}}z^{*} = \frac{\left(\sum_{q=1}^J\frac{1}{p}\right) \left(\sum_{h=r}^K\frac{1}{h}\right)}{t_i}z^{*} = \frac{\sum_{h=r}^K\frac{1}{h}}{t_iK\left(\sum_{p=1}^I\frac{1}{p}\right)}, \quad\forall(i,r)\in\mathcal{U}.
\]

Then, the weights of experts is given by:
\[
W_{i}^{\mathcal{Q}} = \sum\limits_{r=1}^{K} \frac{\sum_{h=r}^K\frac{1}{h}}{t_iK \left(\sum_{p=1}^I\frac{1}{p}\right)} = \frac{1}{t_iK \left(\sum_{p=1}^I\frac{1}{p}\right)} \left(\sum_{r=1}^K\sum_{h=r}^K\frac{1}{h}\right) = \frac1{t_i \left(\sum_{p=1}^I\frac{1}{p}\right)}, \quad \forall i \in \mathcal{U}.
\]

Let $w_{ir}^{*} = W_{i}^{\mathcal{Q}} u_{ir}$, for every $(i,r) \in \mathcal{U}$, we have
\[
    \frac{\sum_{h=r}^K\frac{1}{h}}{t_iK\left(\sum_{p=1}^I\frac{1}{p}\right)} = \frac{u_{ir}}{t_i\left(\sum_{p=1}^I\frac{1}{p}\right)} \Leftrightarrow u_{ir} = \frac{1}{K}\sum_{h=r}^K\frac{1}{h}.
\]

It is evident that the value of $u_{ir}$ is solely determined by the rankings of alternatives, which is a rank-based net utility. It follows that the optimal weights of alternatives under experts, calculated by multiplying the rank-based net utility by the corresponding expert weights, are additive.
\qed
\end{proposition-proof-02}

\section{Proofs for Section \ref{sec-3}}

\newproof{theorem-proof-03}{Proof of Theorem \ref{theorem-03}}
\begin{theorem-proof-03}
Given expert $i$ and attribute $j$, the maximum entropy utility and minimum cross-entropy utility problems have the same partial preference constraints as shown in Equation (\ref{eq-15-2})-(\ref{eq-15-7}), with respective objectives of
\[
\max\limits_{{\mathbf{U}_{ij}}(\cdot)} - \sum\limits_{r=1}^{K_{ij}}U_{ij}(r) \ln U_{ij}(r)
\]
and
\[
    \min \limits_{{\mathbf{U}_{ij}}(\cdot)} \sum\limits_{r=1}^{{{K}_{ij}}}{{{U}_{ij}}(r)\ln \frac{{{U}_{ij}}(r)}{\Bar{V}(r)}}.
\]

Without loss of generality, let $d_{e} \left(U_{ij}(r), U_{ij}(r+1)\right), e\in [E]$ denote the partial preference constraints. Then, using the method of Lagrange multipliers, we have
\[
L^{EU} = - \sum\limits_{r=1}^{K_{ij}}U_{ij}(r) \ln U_{ij}(r) - \lambda_{0} \left( \sum\limits_{r=1}^{{{K}_{ij}}}{{U}_{ij}}(r) - 1 \right) - \sum\limits_{e \in [E]}\lambda_{e} d_{e} \left(U_{ij}(r), U_{ij}(r+1)\right),
\]
and
\[
 L^{CEU} = \sum\limits_{r=1}^{{{K}_{ij}}}{{{U}_{ij}}(r)\ln \frac{{{U}_{ij}}(r)}{\Bar{V}(r)}} + \lambda_{0} \left(\sum\limits_{r=1}^{{{K}_{ij}}}{{U}_{ij}}(r) - 1\right) + \sum\limits_{e \in [E]}\lambda_{e} d_{e} (U_{ij}(r), U_{ij}(r+1)),
\]
where $L^{EU}$ and $L^{CEU}$ denote the Lagrange function of maximum entropy utility and minimum cross-entropy utility, respectively. By taking a partial derivative, we have
\[
    \frac{\partial L^{EU} }{\partial U_{ij}(r)} = 0 \Leftrightarrow U_{ij}^{EU}(r) = e ^ {b(r)}, \quad \forall r \in [K_{ij}],
\]
and
\[
    \frac{\partial L^{CEU}}{\partial U_{ij}(r)} = 0 \Leftrightarrow U_{ij}^{CEU}(r) = \Bar{V}(r) e ^ {b(r)}, \quad \forall r \in [K_{ij}],
\]
where $b(r)$ is the coefficient of partial preference constraints with ranking $r$, which is identical for both maximum entropy utility and minimum cross-entropy utility problems. Under the constraint of normalization, it is evident that maximizing entropy utility in GOPA of discrete prospects with partial preference context is a special case of minimizing cross-entropy utility when the target utility form is uniformly weighted.  
\qed
\end{theorem-proof-03}

\newproof{proposition-proof-03}{Proof of Proposition \ref{proposition-03}}
\begin{proposition-proof-03}
    Given expert $i$ and attribute $j$, by the properties of entropy utility, we can rewrite Equations (\ref{eq-16-2})-(\ref{eq-16-4}) by the expression of indicator function over certain intervals:
    \[
    \begin{aligned}
        & \int_{0}^{K_{ij}} \left(I_{r}(x) - \alpha_{ij{e}_{1}} I_{r-1}(x)\right) u_{ij}(x) \mathrm{d} x = 0, && \quad \forall e_{1} \in [E_{1}], \\
        & \int_{0}^{K_{ij}} \left(I_{r}(x) - I_{r-1}(x) \right) u_{ij}(x) \mathrm{d} x = \beta_{ij{e}_{2}}, && \quad \forall e_{2} \in [E_{2}],\\
        & \int_{0}^{K_{ij}} I_{r}(x) u_{ij}(x) \mathrm{d} x = \gamma_{ij{e}_{3}}, && \quad \forall e_{3} \in [E_{3}],\\
    \end{aligned}
    \]
    where $I_{r}(x)$ are indicator functions with the domain of $[0, r]$.

    Without loss of generality, we can unify the above equations into
    \[
    \int_{0}^{K_{ij}} \zeta_{e}(x) u_{ij}(x) \mathrm{d} x = \xi_{e}, \quad \forall e \in [E].
    \]
    
     It is notable that $\zeta_{e}(x)$ is a step function with jumps at $r_{e}, \forall e\in [E]$.

    Let $F = u_{ij}(x) \ln \big( u_{ij}(x) / v(x) \big), G_{0} = u_{ij}(x), \text{ and } G_{e} = \zeta_{e}(x) u_{ij}(x), \forall e \in [E]$. Construct the auxiliary functional:
    \[
    \begin{aligned}
        J^{*} & = \int_{0}^{K_{ij}} F + \lambda_{0} G_{0} + \sum\limits_{e \in [E]} \lambda_{e} G_{e}  \mathrm{d} x \\
        & = \int_{0}^{K_{ij}} u_{ij}(x) \ln \left(\frac{u_{ij}(x)}{v(x)}\right) + \lambda_{0}u_{ij}(x) + \sum\limits_{e \in [E]} \lambda_{e} \zeta_{e}(x) u_{ij}(x) \mathrm{d} x ,
    \end{aligned}
    \]
    whose Euler equation is
    \[
    \ln\left(\frac{u_{ij}(x)}{v(x)}\right) + 1 + \lambda_{0} + \sum\limits_{e \in [E]} \lambda_{e} \zeta_{e}(x) = 0.
    \]
    
    Rearranging the Euler equation reveals the optimal solution to the minimum cross-entropy utility:
    \begin{equation}
        u^{*}_{ij}(x) = v(x) e^{-1 - \lambda_{0} - \sum\nolimits_{e \in [E]} \lambda_{e} \zeta_{e}(x)},
        \label{eq-sub43-03}
    \end{equation}
    where $\sum\nolimits_{e \in [E]} \lambda_{e} \zeta_{e}(x)$  represents the aggregation of indicator functions, ensuring that $e^{-1 - \lambda_{0} - \sum\nolimits_{e \in [E]} \lambda_{e} \zeta_{e}(x)}$ forms a unique step function over the interval $[0, K_{ij}]$. Thus, the optimal solution to the minimum cross-entropy utility can be expressed as the product of the target utility density function and a unique step function with breakpoints at $r_{e}, \forall e\in [E]$ corresponding to the partial preference information. And the unique step function is also the optimal solution of maximum entropy utility when the target utility density function is uniform (i.e., risk-neutral utility function). \qed
\end{proposition-proof-03}

\newproof{theorem-proof-04}{Proof of Theorem \ref{theorem-04}}
\begin{theorem-proof-04}
    Taking the logarithm of the optimal solution in Proposition \ref{proposition-03} (i.e., Equation (\ref{eq-sub43-03})) yields:
    \[
    \ln (u^{*}_{ij}(x)) = \ln(v(x)) -1 - \lambda_{0} - \sum\nolimits_{e \in [E]} \lambda_{e} \zeta_{e}(x) .
    \]
    
    By Arrow-Pratt’s definition of risk preference function, we have the risk preference function of minimum cross-entropy utility:
    \[
    \eta_{ij}(x) = -\frac{\mathrm{d}}{\mathrm{d}x} \ln (u^{*}_{ij}(x)) = -\frac{\mathrm{d}}{\mathrm{d}x} \ln(v(x)),
    \]
    which is aligned with the risk preference of the target utility density function over each segment. 
    \qed
\end{theorem-proof-04}

\newproof{corollary-proof-03}{Proof of Corollary \ref{corollary-03}}
\begin{corollary-proof-03}
    This follows directly from Proposition \ref{proposition-03}.
    \qed
\end{corollary-proof-03}

\newproof{theorem-proof-05}{Proof of Theorem \ref{theorem-05}}
\begin{theorem-proof-05}
    This follows the symmetric argument as the proof of Theorem \ref{theorem-02}. \qed
\end{theorem-proof-05}

\newproof{corollary-proof-04}{Proof of Corollary \ref{corollary-04}}
\begin{corollary-proof-04}
    Given the case with no missing and duplicate rankings (i.e., $c_{ijr} = 1, \forall (i,j,r) \in \mathcal{U}$), by Equation (\ref{eq-18}) in Theorem \ref{theorem-05}, we have 
    \[
    z^{*}=1 \Bigg/ \left(\sum_{i=1}^I\sum_{j=1}^J\sum_{r=1}^{K_{ij}}\frac{ K_{ij}U^{*}_{ijr}}{t_is_{ij}} \right) = 1 \Bigg/ \left(\sum_{i=1}^I\sum_{j=1}^J\frac{ K_{ij}\sum_{r=1}^{K_{ij}}U^{*}_{ijr}}{t_is_{ij}} \right) = 1 \Bigg/ \left(\sum_{i=1}^I\sum_{j=1}^J\frac{K_{ij}}{t_is_{ij}} \right),
    \]
    where the last equation is from the fact that $\sum_{r=1}^{K_{ij}}U^{*}_{ijr} = 1, \forall(i,j) \in \mathcal{U}$. Hence, $z$ is only related to the number of alternatives and the ranking of experts and attributes. It follows that regardless of the specific global utility structure $v \in \Gamma$, the optimal objective value $z^{*}$ remains unchanged.
    \qed
\end{corollary-proof-04}

\newproof{corollary-proof-05}{Proof of Corollary \ref{corollary-05}}
\begin{corollary-proof-05}
    When there are no missing and duplicate rankings (i.e., $c_{ijr} = 1, \forall (i,j,r) \in \mathcal{U}$), by Corollary \ref{corollary-04}, the optimal value $z^{*}$ is invariant and independent of $U^*_{ijr}$. Then, by Equation (\ref{eq-19}) in Theorem \ref{theorem-05}, we have
    \[
    w_{ij}^{*} = \sum_{r=1}^{K_{ij}} w^{*}_{ijr} = \sum_{r=1}^{K_{ij}} \frac{K_{ij}U^{*}_{ijr}z^{*}}{t_i s_{ij}} = \frac{K_{ij}z^{*} \sum_{r=1}^{K_{ij}}U^{*}_{ijr}}{t_i s_{ij}} = \frac{K_{ij}z^{*}}{t_i s_{ij}}, \quad \forall (i,j) \in \mathcal{U},
    \]
    where $w_{ij}^{*}$ denotes the weight of the attribute $j$ under the preference of the expert $i$. It follows that $w_{ij}^{*}$ is independent of $U^{*}_{ijr}$, which contains the risk preference of experts. The weights of experts and attributes, given by $W^{\mathcal{Q}}_i = \sum_{i=1}^{I}w_{ij}^{*}$ and $W^{\mathcal{N}}_j = \sum_{j=1}^{J}w_{ij}^{*}$, are also independent of $U^{*}_{ijr}$. This indicates that the weights assigned to experts and attributes are independent of the experts' risk preferences toward the attributes.
    \qed
\end{corollary-proof-05}

\newproof{corollary-proof-06}{Proof of Corollary \ref{corollary-06}}
\begin{corollary-proof-06}
    According to Corollary \ref{corollary-05}, $w^{\mathcal{N}}_{j} \text{ and } w_{ij}^{*}$ are independent of the experts' risk preferences toward attributes. Consequently, by Definition \ref{definition-04}, PSD for attributes is also independent of the risk preferences. 
    \qed
\end{corollary-proof-06}

\newproof{lemma-proof-02}{Proof of Lemma \ref{lemma-02}}
\begin{lemma-proof-02}
    The first and second equivalences follow Proposition 2 in \cite{O.R99} and Proposition 2.2 in \cite{D.R03}, respectively.
\qed
\end{lemma-proof-02}

\end{document}